\documentclass{article}

\usepackage{amsmath,amsbsy,amsfonts,amssymb}
\usepackage[francais]{babel}
 \newtheorem{thm}{Th\'eor\`eme}
 
 \newtheorem{cor}[thm]{Corollaire}
   
 \newtheorem{lem}[thm]{Lemme}
  
 \newtheorem{prop}[thm]{Proposition}
   
 \newtheorem{rmq}[thm]{Remarque}
 
 \newtheorem{defi}[thm]{D\'efinition}

\numberwithin{equation}{subsection}

\newcommand{\red}[1]{\rm Red \,}

\title {Paquets stables des s\'eries discr\`etes accessibles par endoscopie tordue; leur param\`etre de Langlands}
\date{}
\author{C. M{\oe}glin}

\begin{document}
\maketitle
\author
\section{Introduction}
Suite aux travaux de J. Arthur en particulier \cite{book}, on sait que la classification de Langlands des repr\'esentations temp\'er\'ees des groupes classiques s'obtient en utilisant l'endoscopie tordue. Ceci est conforme \`a la philosophie de Langlands, l'endoscopie tordue se comprend via la fonctorialit\'e et elle doit donc respecter les paquets stables de repr\'esentations temp\'er\'ees.  Comme l'un des c\^ot\'es est le groupe $GL(n)$ pour lequel la classification de Langlands est connue, l'autre c\^ot\'e s'exprime avec cette classification et on voit avec de l'alg\`ebre lin\'eaire \'el\'ementaire que les groupes endoscopiques tordues se s\'eparent suivant qu'ils sont orthogonaux ou  symplectiques et c'est ainsi que \cite{book} obtient la classification de Langlands pour les groupes symplectiques et orthogonaux (du moins le param\`etre du paquet et il faut toute l'endoscopie d\'evelopp\'ee en \cite{book} pour avoir tout le param\`etre); les groupes unitaires se distinguent aussi (avec les fonctions $L$ d'Asai) (les m\'ethodes de \cite{book} sont reprises dans ce cadre par \cite{mok}); les groupes $GSpin$ s'obtiennent aussi de cette fa\c{c}on (ceci a \'et\'e expliqu\'e en \cite{gsp4}). Il se peut que d'autres situations avec des changements de base s'obtiennent aussi ainsi et c'est pour cela que l'on veut pr\'esenter ici les arguments les plus simples et les plus g\'en\'eraux pour arriver \`a ce r\'esultat.

La preuve de ce r\'esultat, pr\'esent\'ee ici, est  essentiellement locale; la seule partie globale est le fait qu'une identit\'e de caract\`eres entre repr\'esentations elliptiques qui est satisfaite pour les fonctions dont les int\'egrales orbitales sont nulles sur les \'el\'ements non elliptiques se prolonge en une identit\'e de caract\`ere (sans restriction de support), r\'esultat d\'emontr\'e dans le cas non tordu en \cite{selecta}. Ce r\'esultat n\'ecessite la formule des traces simple et sa stabilisation: il faut en effet d\'emontrer qu'une telle identit\'e de caract\`ere se prolonge en une identit\'e de caract\`ere quitte \`a y ajouter des repr\'esentations induites \`a partir de sous-groupes paraboliques propres, c'est la m\'ethode de \cite{selecta}. Et c'est \`a cet endroit qu'on utilise la formule des traces. Ensuite des m\'ethodes locales permettent d'enlever ces repr\'esentations suppl\'ementaires. Cette m\'ethode a \'et\'e expliqu\'ee en \cite{waldsgl} (reprise en \cite{pacific}) et on l'\'ecrira en toute g\'en\'eralit\'e pour les besoins de la stabilisation de la formule des traces tordues dans \cite{avenir}. Un tel r\'esultat est un pr\'eliminaire \`a toute classification \`a la Langlands des repr\'esentations temp\'er\'ees (\cite{book} et \cite{mok} en d\'ependent); et le but de cet article est d'exploiter compl\`etement un tel r\'esultat. 

On ne peut pas dire que les m\'ethodes ci-dessous soient r\'eellement diff\'erentes de celles de \cite{book} mais elles sont plus ax\'ees sur la th\'eorie des repr\'esentations et on obtient donc des cons\'equences plus fines en th\'eorie des repr\'esentations.

On a clairement distingu\'e ce qui peut se d\'emontrer avec le r\'esultat de prolongement de l'identit\'e de caract\`eres expliqu\'e ci-dessus et des propri\'et\'es \`a peu pr\`es \'el\'ementaires de th\'eorie de repr\'esentation, de ce qui n\'ecessite l'introduction du $L$-groupe et qui est surtout de l'alg\`ebre lin\'eaire \'el\'ementaire. Le lien entre les deux est la doubling method (cf. ci-dessous) qui elle n'est pas \'el\'ementaire. Plus pr\'ecis\'ement, la classification cherch\'ee peut se d\'ecomposer en deux \'etapes.

Dans la premi\`ere partie, on d\'emontre qu'\'etant donn\'e $\pi$ une s\'erie discr\`ete d'un groupe classique, la projection d'un pseudo coefficient (cuspidal) de $\pi$ sur la partie stable (cf. \ref{pseudocoefficient}) d\'etermine une unique repr\'esentation virtuelles du groupe, $\pi_{st}$, combinaison lin\'eaire uniquement form\'ee de s\'eries discr\`etes, stable par construction et que le transfert tordu de cette repr\'esentation virtuelle est une repr\'esentation irr\'eductible (\`a un scalaire pr\`es), $\pi^{GL}$ du groupe $GL(n)$ de la situation; c'est le mot irr\'eductible ici qui est important. Le support cuspidal de cette repr\'esentation de $GL(n)$ est l'invariant des paquets stables (cf \ref{supportcuspidaletendu}); c'est assez agr\'eable car le support cuspidal se contr\^ole \'evidemment tr\`es bien par induction et restriction et qu'il est d\'etermin\'e par des propri\'et\'es de r\'eductibilit\'e d'induites. On en revient \`a notre d\'efinition des blocs de Jordan, que l'on modifie l\'eg\`erement \`a cet endroit par commodit\'e mais on d\'emontre finalement (en \ref{unicitedefinition}) que l'on n'a en fait pas chang\'e la d\'efinition originelle de \cite{europe} et \cite{jams}

On montre aussi facilement que si $\pi$ et $\pi'$ sont des s\'eries discr\`etes du groupe classique telle que $\pi_{st}$ n'est pas orthogonal \`a $\pi'_{st}$ alors $\pi_{st}$ est proportionnel \`a $\pi'_{st}$ et \`a un scalaire les repr\'esentations $\pi^{GL}$ et $\pi^{'GL}$ co\"{\i}ncide. R\'eciproquement si $\pi^{GL}$ est une repr\'esentation temp\'er\'ee (elliptique quand on prend en compte l'automorphisme de la situation) du groupe $GL(n)$ il existe une unique donn\'ee endoscopique ellitpique tel que $\pi^{GL}$ soit un transfert d'une des repr\'esentations $\pi_{st}$ (\`a un scalaire pr\`es) obtenue par projection d'une s\'erie discr\`ete du groupe endoscopique sous-jacent (comme pr\'ec\'edemment).

On obtient ainsi une classification des combinaisons lin\'eaires stables de s\'eries discr\`etes \`a l'aide des groupes $GL(n)$. En utilisant la correspondance de Langlands pour les groupes $GL(n)$, on a donc associ\'e \`a toute s\'erie discr\`ete d'un groupe classique un morphisme de $W_{F}\times SL(2,{\mathbb C})$ dans le $L$-groupe de $GL(n)$.

Dans la deuxi\`eme partie de l'article on montre que ce morphisme est \`a valeurs dans le $L$-groupe du groupe classique. Et on a donc ainsi une classification \`a la Langlands des paquets stables.

Le lien entre ces deux parties est la doubling method dont Pia\-te\-tski-Sha\-piro a \'et\'e avec en particulier S. Rallis l'un des concepteurs (r\'edig\'e par J. Cogdell dans \cite{GPSR}.   La doubling method s'applique \`a la th\'eorie des repr\'esentations et permet de relier les points de r\'eductibilit\'e de certaines induites \`a des propri\'et\'es de fonctions $L$ d\'efinies par Shahidi. 
Et ces fonctions $L$ de Shahidi  sont aussi les fonctions $L$ du groupe $W_{F}$ gr\^ace aux r\'esultats d'Henniart (\cite{henniartasai}); la preuve de \cite{henniartasai} est en grande partie globale, elle utilise l'\'equation fonctionnelle. Et du c\^ot\'e galoisien la doubling method donne une factorisation des fonctions $L$ qui donne exactement la dichotomie, orthogonal/symplectique, Asai/Asai tordu. 

Les p\^oles des fonctions $L$ gouvernent la factorisation du morphisme de Langlands associ\'e \`a une repr\'esentation d'un groupe g\'en\'eral lin\'eaire dans le groupe dual des groupes consid\'er\'es ici; c'est \'evident quand ces fonctions sont des fonctions $Sym^2$ ou $\wedge^2$ et c'est aussi le cas par exemple pour les fonctions $L$ d'Asai (comme cela introduit des signes moins concrets, on refait en d\'etail les calculs).

Et c'est ce qui fait le lien entre la classification purement en terme de repr\'esentations et la construction du morphisme de Langlands associ\'e \`a un paquet stable de s\'eries discr\`etres. La premi\`ere partie se termine en \ref{conclusion}, la deuxi\`eme partie d'alg\`ebre lin\'eaire (fortement ax\'e sur les groupes unitaires) fait l'objet du paragraphe 5 et le lien fourni par la doubling method est le paragraphe 6.

\

On n'a pas calcul\'e les coefficients dans les formules de transfert alors que \cite{selecta} et\cite{book} trace la voie pour le faire et le fait pour les groupes orthogonaux et symplectiques; White (\cite{white}) et Mok \cite{mok} ont \'etendu ces travaux au cas des groupes unitaires. La m\'ethode n'est pas locale, il faut mettre la situation locale dans une situation globale o\`u on conna\^{\i}t les r\'esultats en toutes les places sauf celle qui nous int\'eresse. En particulier, il faut conna\^{\i}tre (dans le point de vue de \cite{book}) les r\'esultats aux places archim\'ediennes, r\'esultats en cours de d\'emonstration par Mezo; pour les groupes unitaires  les r\'esultats n\'ecessaires aux places archim\'ediennes sont d\'ej\`a connus gr\^ace aux travaux de L. Clozel \cite{clozel}.

Nous ne faisons pas non plus la classification fine de Langlands de toutes les s\'eries discr\`etes; ceci est uniquement pour ne pas allonger l'article, les m\'ethodes locales que nous avons d\'evelopp\'ees dans \cite{europe}, \cite{jams}, \cite{pacific} et \cite{pourshahidi} (partiellement en collaboration avec M. Tadic) s'appliquent sans probl\`eme. De m\^eme on s'est limit\'e au cas des groupes quasi-d\'eploy\'es ce qui n'est pas n\'ecessaire.

Comme on obtient les repr\'esentations temp\'er\'ees comme module de Jacquet de certaines s\'eries dis\-cr\`etes, la classification des s\'eries dis\-cr\`etes est suffisante pour avoir une classification des paquets stables de repr\'esen\-tations temp\'er\'ees; toutefois le $R$-groupe joue un r\^ole important car on ne peut travailler avec les s\'eries discr\`etes sans travailler avec les repr\'esentations elliptiques (qui sont des repr\'esentations temp\'er\'ees, virtuelles d\'efines avec le $R$-groupe \cite{arthurelliptic}); ces $R$-groupes ont \'et\'e calcul\'es en particulier en \cite{herb} et \cite{goldberg}.

Ce fut un r\'eel honneur pour moi d'\^etre invit\'ee \`a participer au colloque ''Legacy of I. I. Piatetski-Shapiro'' ayant eu lieu \`a l'Universit\'e de Yale du 6 au 12 Avril 2012. Je remercie tr\`es chaleureusement les organisateurs de cette conf\'erence et d\'edie ce texte \`a la m\'emoire de Ilya Piateskii-Shapiro.

\tableofcontents
\section{Notations et propri\'et\'es g\'en\'erales}
\subsection{Les groupes\label{groupes}}
Ici $F$ est un corps p-adique et $E$ est une extension de $F$ soit \'egale \`a $F$ soit une extension quadratique de $F$. On note $\tilde{GL}$ l'un des groupes alg\'ebriques suivant: $GL(n,E)$, $GL(n,F)\times F^*$ et $\theta$ l'automorphisme ext\'erieur $$\theta(g,\lambda)
:=( \, ^t\overline{g}^{-1},\lambda det\, g),$$
o\`u il n'y a pas de deuxi\`eme facteur dans le cas $GL(n,E)$ et o\`u $\overline{g}$ est le conjugu\'e de $g$ si $E\neq F$.

Dans la suite de cet article on note $\underline{G}$ une donn\'ee endoscopique elliptique de la composante non connexe du produit semi-direct de $\tilde{G}$ avec $\{1,\theta\}$. Ce qui nous int\'eresse est le groupe sous-jacent \`a cette donn\'ee que l'on note $G_{n}$; on pourra faire dispara\^{\i}ntre l'indice $n$ s'il n'y a pas de confusion possible. Les groupes $G$ que l'on obtient ainsi, sont quasi-d\'eploy\'es et leur ensemble contient la liste: $Sp(n-1,F)$ pour $n$ impair, $SO(n,F)$ et $SO(n+1,F)$ avec $n$ pair (toutes les formes quasi-d\'eploy\'ees) , $GSpin_{n}(F)$ (toutes les formes quasi-d\'eploy\'ees)  et $GSpin_{n+1}(F)$ avec $n$ pair et $U(n,E/F)$. Les groupes de cette liste sont appel\'es, les groupes endoscopiques simples. Dans presque tous les cas, le groupe d\'etermine la donn\'ee endoscopique mais ce n'est pas vrai ni pour $Sp(n-1,F)$ ni pour $U(n,E/F)$; le cas de $U(n,E/F)$ est bien connu puisqu'il y a le cas de ce que l'on appelle improprement le changement de base stable et instable. Le cas de $Sp(n-1,F)$ se trouve expliqu\'e dans \cite{book}; il faut tenir compte de l'action du groupe de Galois qui fournit, en plus,  un caract\`ere quadratique non n\'ecessairement trivial.

A cette liste s'ajoute  le groupe endoscopique principal pour $\tilde{GL}.\theta$ quand $\tilde{GL}=GL(n,F)\times F^*$ avec $n$ impair; c'est $Sp(n-1,F)\times F^*$. Et s'ajoute les groupes endoscopiques qui sont des produits de  groupes dans la liste pr\'ec\'edente.

Une description pr\'ecise que nous utiliserons se trouve dans \cite{inventaire} dans le cas o\`u $\tilde{GL}$ n'a pas le facteur $F^*$ et ce dernier cas est d\'ecrit dans \cite{gsp4}.

On notera \`a quelques endroits $GSpin_{2n}^*(F)$ le groupe non connexe qui est un rev\^etement du groupe de similitudes orthogonales sur un espace de dimension 2n. C'est naturellement un rev\^etement de $O(2n,F)$.

\subsection{Les repr\'esentations}
On s'int\'eresse aux repr\'esentations de $\tilde{GL}$ qui sont stables sous $\theta$; pour simplifier un peu les notations, on fixe un caract\`ere $\nu$ de $F^*$, unitaire, et pour toute repr\'esentation $\pi$ de $GL(n,E)$, on note $\tilde{\pi}$ sa contragr\'ediente si $E=F$ et son dual hermitien si $E\neq F$. Soit $\pi$ une repr\'esentation de $\tilde{GL}$ telle que $\theta.\pi \simeq \pi$. On suppose que la restriction de $\pi$ \`a $F^*$ (quand ce facteur existe) est le caract\`ere $\nu$. La $\theta$ invariance se traduit donc par:
$$
\tilde{\pi}\otimes \nu \simeq \pi.\eqno(1)
$$

La notion de s\'eries discr\`etes est bien connue, la caract\'erisation utilis\'ee ici est le fait que les exposants sont dans une chambre de Weyl obtuse, ouverte positive. Les repr\'esentations temp\'er\'ees sont celles dont les exposants sont dans une chambre de Weyl obtuse positive ferm\'ee. Les repr\'esentations elliptiques ont \'et\'e d\'efinies par Harish-Chandra mais on utilise la variante de \cite{arthurelliptic} (premiers paragraphes de cet article). On a besoin de ces notions dans le cas tordu; cela est fait en g\'en\'eral dans \cite{ftl} 2.12 mais ici nous n'en avons besoin essentiellement pour $\tilde{GL}.\theta$. La description se fait alors en termes tr\`es simples: il n'y a pas de diff\'erence entre $\theta$-elliptique et $\theta$-discr\`ete, ce sont des repr\'esentations irr\'eductibles $\theta$-invariantes induites de repr\'esentations de Steinberg toutes non isomorphes. On dira indiff\'eremment, $\theta$-elliptique ou $\theta$-discr\`ete en pr\'ef\'erant la premi\`ere terminologie.

Tr\`es accessoirement, on aura besoin de la notion de repr\'esentation elliptique pour les groupes $O(2n,F)$ et $GSpin^*_{2n}(F)$; on l'utilisera dans un cas tr\`es simple o\`u on induit \`a partir du parabolique maximal de Levi isomorphe \`a $GL(n,F)$ ou $GL(n,F)\times F^*$ une repr\'esentation $\theta$ invariante. Quelque soit la parit\'e de $n$, l'induite au groupe non connexe se coupe en deux et la diff\'erence des deux repr\'esentations est une repr\'esentation elliptique pour la composante non neutre de ces groupes.

R. Herb (\cite{herb}) a donn\'e une caract\'erisation des repr\'esentations elliptiques dans le cas des groupes que nous g\'en\'eraliserons au cadre tordu dans le cadre de la stabilisation de la formule des traces tordu.
\subsection{Les pseudo-coefficients\label{pseudocoefficient}}
On appelle fonctions cuspidales sur $G$ une fonction lisse dont les int\'egrales orbitales pour les \'el\'ements semi-simples r\'eguliers hyperboliques (c'est-\`a-dire inclus dans un sous-groupe parabolique propre de $G$) sont nuls. Et
on note $I_{cusp}(G)$ l'espace des fonctions cuspidales sur $G$ modulo celles dont toutes les int\'egrales orbitales sur des \'el\'ements semi-simples r\'eguliers sont nulles. On adopte la m\^eme notation pour $\tilde{G}$ \'etant entendu que l'on consid\`ere l\`a des fonctions sur la composante $\tilde{G}.\theta$; la notation de parabolique est alors remplac\'ee par la notion d'espace parabolique ce qui pompeusement d\'esigne les sous-groupes paraboliques $\theta$-stable de $G$.

En \cite{selecta} \'etendu en \cite{ftl}, paragraphe 7, pour inclure le cas tordu, il est montr\'e  que $I_{cusp}^G$ et $I_{cusp}^{\tilde{GL}}$ s'interpr\`etent comme l'espace des pseudo-coefficients des repr\'esentations elliptiques de $G$ ou $\tilde{G}$ (pour le cas tordu, \cite{ftl} 7.2 (1)). On peut stabiliser chacun de ces espaces (\cite{selecta} pour le cas non tordu et \cite{waldsstabilisation} 4.11 (i) de la proposition  pour le cas tordu), ce qui permet de d\'efinir $I_{cusp, st}^G$ la partie stable et de le faire aussi pour tous les groupes endoscopiques de $G$; quand on r\'ealise $G$ comme le groupe dans une donn\'ee endoscopique $\underline{G}$ de $\tilde{GL}.\theta$, il faut aussi tenir compte des automorphismes de cette donn\'ee. On \'ecrit alors $I_{cusp,st}^{\underline{G}, Aut}$ pour l'image du transfert de $I_{cusp}^{\tilde{GL}}$. Pour $\tilde{GL}$, on a:
$$
I_{cusp}^{\tilde{G}}=\oplus_{\underline{G}}I_{cusp,st}^{\underline{G}, Aut},\eqno (1)
$$
ici $\underline{G}$ parcourt toutes les donn\'ees endoscopiques elliptiques de $\tilde{GL}.\theta$ et  l'\'egalit\'e est donn\'ee par la somme des transferts. 

Pour $G$ fix\'e comme ci-dessus, on \'ecrit la stabilisation de \cite{selecta}
$$
I_{cusp}^{G}=\oplus_{\underline{H}}I_{cusp,st}^{\underline{H}, Aut},\eqno(2)
$$
qui est l'analogue non tordu de (1) et o\`u $\underline{H}$ parcourt l'ensemble des donn\'ees endoscopiques elliptiques de $G$ y compris $G$ lui-m\^eme.

On dit qu'une repr\'esentation virtuelle combinaison lin\'eaire de repr\'e\-senta\-tions elliptiques, $\pi$ de $G$ est stable sur les elliptiques si la combinaison lin\'eaire des pseudo-coefficients dans $I_{cusp}^G$ est en fait un \'el\'ement de $I_{cusp,st}^{G}$; on identifie $\pi$ \`a cette combinaison lin\'eaire de pseudo-coefficients et donc \`a un \'el\'ement de $I_{cusp}^{G}$. Cela est \'equivalent \`a dire que la distribution $f\in I_{cusp}^{G}\mapsto tr\, \pi(f)$ est stable. 

Soit $\pi$ un \'el\'ement de $I_{cusp}^{G}$ vu comme une repr\'esentation virutelle; en utilisant (2),
on d\'etermine pour toute donn\'ee endoscopique elliptique $\underline{H}$ de $G$, une repr\'esentation virtuelle $\pi_{\underline{H},st}$ telles que pour tout $f\in I_{cusp}^{G}$, on ait l'\'egalit\'e:
$$
tr\, \pi(f)=\sum_{\underline{H}}tr\, \pi_{\underline{H},st}(f^{\underline{H}}),\eqno(3)
$$
o\`u $f^{\underline{H}}$ est un transfert de $f$.

Un r\'esultat cl\'e de \cite{selecta} est de montrer que si $\pi$ est dans $I_{cusp,st}^G$, alors la distribution $tr\, \pi$ est stable au sens usuel et que, sous les hypoth\`eses de l'\'egalit\'e (3), alors cette \'egalit\'e se prolonge \`a toute fonction lisse $f$.

La g\'en\'eralisation de ces propri\'et\'es au cas tordu est un pr\'eliminaire \`a tous les travaux qui \'etablissent les propri\'et\'es locales des repr\'esentations temp\'er\'ees (\cite{book},\cite{mok}). G\'en\'eraliser \cite{selecta} avec les m\^emes m\'ethodes est s\^urement possible mais il faudrait r\'ecrire des centaines de pages \'ecrites par Arthur dans le cas non tordu. En \cite{waldsgl} et \cite{pacific} on a obtenu ces r\'esultats avec une autre m\'ethode. Cela ne couvre pas tous les cas que l'on a en vu ici mais nous avons maintenant v\'erifi\'e (\cite{avenir}) que cette autre m\'ethode s'\'etend en toute g\'en\'eralit\'e au cas tordu; on admet  ici ce r\'esultat (qui sera disponible sous forme de pr\'epublication tr\`es prochainement). Donc on admet que si $\pi^{\tilde{GL}}\in I_{cusp}^{\tilde{GL}}$ est \'ecrit suivant la somme directe (1) en $\pi^{\tilde{GL}}=\oplus_{\underline{G}}\pi_{\underline{G},st}$, on a pour toute fonction lisse $f$ sur $\tilde{GL}.\theta$ l'\'egalit\'e des traces
$$
tr\, \pi^{\tilde{GL}}(f)=\sum_{\underline{G}}tr\, \pi_{\underline{G},st}(f^{\underline{G}}).
$$
\subsection{S\'eries discr\`etes et paquets stables\label{projectionnonnulle}}
\begin{prop} Soit $\pi$ une s\'erie discr\`ete de $G$ vue comme un \'el\'ement de $I_{cusp}^G$ (via un pseudo coefficient cuspidal). Alors la projection de $\pi$ sur $I_{cusp,st}^{G}$ est non nulle.
\end{prop}
La d\'emonstration de cette proposition est d\'ej\`a dans \cite{pacific} paragraphe 2 et m'a \'et\'e donn\'ee par Waldspurger: on regarde les germes du caract\`ere au voisinage de l'origine; ils se d\'eveloppent par degr\'e d'homog\'en\'eit\'e et le degr\'e formel est l'un de ces termes. Ce terme est stable et sa projection stable est donc non nulle et il y a donc un germe de la projection de  qui est non nul. D'o\`u la non nullit\'e de la projection
\begin{rmq} Dans cet article on montre que la proposition est stric\-tement sp\'e\-cifique aux s\'eries discr\`etes, c'est-\`a-dire que la projection sur $I_{cusp,st}^G$ d'une repr\'e\-sentation elliptique combinaison lin\'eaire de repr\'e\-sentations temp\'er\'ees dont aucune n'est une s\'erie discr\`ete est nulle. Une d\'emonstration a priori de ce r\'esultat simlifierait les d\'emonstrations mais je ne sais pas si une telle propri\'et\'e est vraie en dehors des groupes consid\'er\'es ici.
\end{rmq}
On appelle paquet stable contenant $\pi$ la repr\'esentation virtuelle, somme de repr\'esentations elliptiques, stable et dont l'image dans $I_{cusp,st}^G$ est \'egale \`a la projection du pseudo coefficient de $\pi$.
\subsection{Notations pour les modules de Jacquet\label{jac}}
Soit $\pi$ une repr\'esentation  de $G_{n}$ et soit $d$ un entier; on suppose que $G_{n}$ contient un sous-groupe parabolique de Levi isomorphe \`a $GL(d,E)\times G_{n-2d}$ et on fixe $\rho$ une repr\'esentation cuspidale unitaire irr\'eductible de $GL(d,E)$. On fixe aussi $x$ un nombre r\'eel et on note $Jac_{\rho\vert\,\vert^x}\pi $ la repr\'esentation virtuelle dans le groupe de Grothendieck de $G_{n-2d}$ tel que le module de Jacquet de $\pi$ pour un parabolique de sous-groupe de Levi isomorphe \`a $GL(d,E)\times G_{n-2d}$ soit la somme de $\rho\vert\,\vert^x\otimes Jac_{\rho\vert\,\vert^x}\pi$ et de repr\'esentation irr\'eductible de la forme $\sigma\otimes \tau$ o\`u $\sigma$ est une repr\'esentation de $GL(d,E)$ non isomorphe \`a $\rho\vert\,\vert^x$.

On adopte une notation analogue pour $\tilde{GL}_{n}$: ici on fixe un espace para\-bolique   de $\tilde{GL}$, c'est-\`a-dire un sous-groupe parabolique, $P$ de $\tilde{GL}$ dont le normalisateur dans $\tilde{GL}.\theta$ est non trivial. On fixe alors un espace de Levi, $M,\tilde{M}$ de cet espace parabolique et on suppose que le sous-groupe de Levi sous-jacent est isomorphe \`a $GL(d,E)\times \tilde{GL}_{n-2d}\times GL(d,E)$; on fixe les notations en fixant $w$ un \'el\'ement de $\tilde{GL}$ qui laisse stable $M$ et conjugue les deux copies de $GL(d,E)$ et on pose $\theta':=w \theta$. Alors $\tilde{M}=M.\theta'$.  Soit $\pi$ une repr\'esentation de $\tilde{GL}_{n}$ prolong\'ee en $\tilde{\pi}$,  \`a $\tilde{GL}_{n}.\theta$; on d\'efinit le module de Jacquet de $\tilde{\pi}$; c'est aussi le module de Jacquet de $\pi$ mais il a une action canonique de $\theta$ et donc de $GL(d,E)\times \tilde{GL}_{n-2d}.\theta_{n-2d} \times GL(d,E)$, l'espace de Levi. On fixe $\rho$ une repr\'esentation de $GL(d,E)$, cuspidal unitaire irr\'eductible et on suppose que $\rho\simeq \tilde{\rho}\otimes \nu$ (cf \ref{groupes} (1)).

On \'ecrit cette repr\'esentation comme la somme d'une repr\'esentation, $\pi_{\rho}$ dont tous les sous-quotients irr\'eductibles comme repr\'esentation de $GL(d,E)\times \tilde{GL}_{n-2d}\times GL(d,E)$ sont de la forme $\rho\vert\,\vert^x \times \sigma\times \rho\vert\,\vert^{-x}$ (o\`u $\sigma$ est queconque) et une autre dont aucun sous-quotient irr\'eductible n'a cette propri\'et\'e. On remarque que cette d\'ecomposition est stable  sous $\theta'$ et est canoniquement une repr\'esentation de $\tilde{M}$. On note $Jac_{\rho\vert\,\vert^x}^{\theta}\pi$ une repr\'esentation $\theta_{n-2d}$ invariante (vue dans le groupe de Grothendieck) telle que $\rho\vert\,\vert^x \otimes Jac_{\rho\vert\,\vert^x}^{\theta}\pi \otimes \rho\vert\,\vert^{-x}$ et $\pi_{\rho}$ ont m\^eme trace sur $\tilde{M}$.
\subsection{Modules de Jacquet et transfert\label{transfertetjac}}
Depuis les travaux pionniers de Shelstad, on sait que le tranfert commute \`a l'induction: quand $\pi$ est une repr\'esentation induite d'une repr\'esentation $\sigma$, il est facile de calculer la trace de $\pi$ en fonction de la trace de $\sigma$ vue comme une repr\'esentation d'un sous-groupe de Levi $M$. En effet, $tr\, \pi(f)=tr\, \sigma(f_{M})$, o\`u $f_{M}$ est la fonction sur $M$ d\'efinie par $f(m)=\int_{K\times N}f(k^{-1}nmk)$ avec les notations usuelles (il faut ajouter un facteur discriminant). On se place dans la situation o\`u $\underline{G}'$ est une donn\'ee endoscopique elliptique de $G$. On fixe $M$ un sous-groupe de Levi de $G$ qui provient par transfert d'un sous-groupe de Levi $\underline{M}'$ de $\underline{G}'$. Pour savoir que l'induction commute au transfert, il est alors suffisant de v\'erifier que si $f'$ est un transfert de $f$ alors $f'_{M'}$ est un transfert de $f_{M}$; un calcul \'el\'ementaire sur les int\'egrales orbitales ram\`enent l'assertion au fait de savoir  que les facteurs de transfert pour les \'el\'ements de $M,\underline{M'}$ sont les m\^emes qu'on les calcule dans $G,\underline{G}'$ ou dans $M,\underline{M}'$ et ceci est fait en \cite{shelstad} 3.4.2 au moins pour les groupes r\'eels.  Le cas des groupes p-adiques et de l'induction dans la situation tordue (il faut alors bien \'evidemment que l'induction se fasse via un parabolique stable) est analogue. Cette propri\'et\'e relative \`a l'induction est utilis\'ee partout dans les travaux sur la formule des travaux et en particulier dans \cite{selecta}.

Ici on pr\'ef\`ere utiliser le fait que le module de Jacquet commute au transfert; cela repose sur la m\^eme propri\'et\'e des facteurs de transfert. En effet, la trace sur les modules de Jacquet s'exprime  en fonction de  la trace sur la repr\'esentation en des points $tm$ o\`u $t$ est un \'el\'ement du centre de $M$ qui contracte suffisamment le radical unipotent du parabolique (cf \cite{mwselecta}) 4.2.1. En loc. cite, on ne consid\'erait que des groupes endoscopiques principaux, les facteurs de transfert sont alors \'egaux \`a $1$ et l'\'egalit\'e est triviale; d'autres cas  tels que le transfert spectral s'obtient \`a partir de l'endoscopie principale par torsion par un caract\`ere (le cas du transfert non principal du groupe unitaire par exemple) s'en d\'eduisent aussi.  Le cas g\'en\'eral vient de l'\'egalit\'e du paragraphe pr\'ec\'edent (la preuve de 6.2.1 de \cite{taibi} \'ecrit le cas des groupes orthogonaux pairs qui est avec $GSpin_{2n}$ l'unique cas non trivial pour ce que l'on fait ici).

En conclusion, avec les notations du paragraphe pr\'ec\'edent, on fixe $\pi^{GL}$ une repr\'esentation de $\tilde{GL}$ et $M$ un espace de Levi de c'est-\`a-dire un sous-groupe de Levi d'un parabolique $\theta$-stable. On \'ecrit en fonction de \ref{pseudocoefficient} (1) la projection de la trace tordue de $\pi^{GL}$ sur chaque $I_{cusp}^{\underline{G},s}$:
$$
tr_{\theta}\, \pi^{GL}=\oplus_{\underline{G}} tr\, \pi_{\underline{G},st}
$$
et on fixe $M,\rho,x$ comme dans \ref{jac}
$$
Jac^{\theta}_{\rho\vert\,\vert^{x}}\pi^{GL}=\oplus_{\underline{G}} Jac_{\rho\vert\,\vert^x}\pi_{\underline{G,st}},
$$
et ceci est une \'egalit\'e de transfert mais dans le terme de droite chaque repr\'esen\-tation est une repr\'esentation virtuelle d'un sous-groupe de Levi d'un groupe endoscopique elliptique de $\tilde{GL}.\theta$ et certaines se regroupent puisque que plusieurs donn\'ees endoscopiques elliptiques peuvent avoir m\^eme donn\'ee endoscopique pour le sous-groupe de Levi se transf\'erant en $M$.

On verra au d\'ebut de la preuve du th\'eor\`eme de \ref{projectionendoscopique} qu'aucun regroupement n'a lieu si $x>1/2$.

\subsection{Quelques propri\'et\'es g\'en\'erales\label{pasdeuxfois}}
\begin{prop}
Soit $\pi$ une s\'erie discr\`ete de $G$; on fixe $\rho$ une repr\'esentation cuspidale unitaire de $GL(d_{\rho},E)$ et $x\in {\mathbb R}-\{0\}$. Alors 
$Jac_{\rho\vert\,\vert^x}\circ Jac_{\rho\vert\,\vert^x} \pi=0.$
\end{prop}
Avant de faire la preuve, remarquons que l'hypoth\`ese $x\neq 0$ pourra \^etre enlev\'ee \`a la fin de l'article mais en ce d\'ebut d'article, on l'utilise dans la  preuve.

On consid\`ere la repr\'esentation virtuelle elliptique stable contenant $\pi$ qui est d\'ecrite en \ref{projectionnonnulle}, on la note $\pi_{st}$. On fixe $t$ maximal avec le fait qu'il existe une repr\'esentations irr\'eductible $\pi'$ intervenant dans ce paquet (soit une s\'erie discr\`ete soit une composante d'une repr\'esentation elliptique) pour laquelle appliqu\'e $t$ fois $Jac_{\rho\vert\,\vert^x}$ donne un r\'esultat non nul. On note $\pi'$ une telle repr\'esentation et il suffit de montrer que $t\leq 1$.
Par r\'eciprocit\'e de Frobenius, on sait qu'il existe une repr\'esentation irr\'eductible $\tau'$ de $G_{n-2td_{\rho}}$ et une inclusion
$$
\pi'\hookrightarrow \rho\vert\,\vert^x\times \cdots \times \rho\vert\,\vert^{x}\times \tau'.\eqno(1)
$$
Par maximalit\'e de $t$, $Jac_{\rho\vert\,\vert^x}\tau'=0$. Ainsi $Jac_{\rho\vert\,\vert^x}$ appliqu\'e $t$ fois \`a $\pi'$ donne un multiple de $\tau'$. On v\'erifie que $\tau'$ n'est sous-quotient d'aucun $Jac_{\rho\vert\,\vert^x}$ appliqu\'e $t$ fois \`a l'une des autres s\'eries discr\`etes intervenant dans le paquet consid\'er\'e; en effet si $\pi''\neq \pi'$ avait cette propri\'et\'e, on aurait d'abord l'existence de $\tau''$ avec une inclusion analogue \`a (1) et aussi, par maximalit\'e de $t$, $Jac_{\rho\vert\,\vert^x}\tau''=0$. Ainsi n\'ecessairement, parce que $x\neq 0$,  $\tau''\simeq \tau'$ et $\pi''$ est aussi un sous-module irr\'eductible de (1). On consid\`ere la repr\'esentation de $GL(td_{\rho},E)$ isomorphe \`a l'induite irr\'eductible de $t$ facteurs $\rho\vert\,\vert^x$ et on la note $\sigma$. On remarque que $\sigma\otimes \tau$ a multiplicit\'e un dans le module de Jacquet de $\pi$ comme quotient irr\'eductible; cela vient \'evidemment du fait que $x\neq 0$ et $Jac_{\rho\vert\,\vert^x}\tau=0$. Par r\'eciprocit\'e de Frobenius, cela assure que (1) a un unique sous-module irr\'eductible. Ainsi $\pi''=\pi'$ et quand on applique $t$ fois $Jac_{\rho\vert\,\vert^x}$ \`a la distribution stable consid\'er\'ee on obtient une repr\'esentation virtuelle non nulle, not\'ee $Jac_{\rho\vert\,\vert^x, t-fois}\pi_{st}$. 

On note $\pi^{GL}$ le transfert de $\pi_{st}$.
Alors $Jac_{\rho\vert\,\vert^x, t-fois}\pi_{st}$
se transfert en $Jac_{\rho\vert\,\vert^x}^{GL}$ appliqu\'e $t$ fois \`a $\pi^{GL}$. Donc ce module de Jacquet doit  lui aussi \^etre non nulle mais ceci est impossible si $t>1$ puisque $\pi^{GL}$ est une combinaison lin\'eaire de repr\'esentations $\theta$-elliptiques irr\'eductibles; une repr\'esentation $\theta$-elliptique irr\'eductible est une induite de repr\'esentation de Steinberg pour un sous-groupe de Levi, toutes les repr\'esentations de Steinberg \'etant in\'equivalentes. D'o\`u la proposition.
\begin{cor}(i) Soit $\pi$ une composante d'une repr\'esentation elliptique intervenant dans un paquet stable et soit  $x\in {\mathbb R}-\{0\}$. Soit $Jac_{\rho\vert\,\vert^x}\pi =0$, soit $x>0$ et $Jac_{\rho\vert\,\vert^{x}}\pi$ est une repr\'esentation irr\'eductible.

(ii) Soient $\pi,\pi'$ deux repr\'esentations in\'equivalentes de $G$ intervenant comme composantes dans un paquet stable de repr\'esentations elliptiques et soient $\rho,x$ comme en (i). Alors soit $Jac_{\rho\vert\,\vert^x}\pi\neq  Jac_{\rho\vert\,\vert^x}\pi'$ soit $Jac_{\rho\vert\,\vert^x}\pi= Jac_{\rho\vert\,\vert^x}\pi'=0$.
\end{cor}
Ce corolaire a \'et\'e d\'emontr\'e dans la preuve ci-dessus puisque $t=1$, avec les notations de cette preuve.
\subsection{Propri\'et\'es de la projection endoscopique\label{projectionendoscopique}}
Soit $\pi^{GL}$ une repr\'esentation $\theta$-elliptique irr\'eductible de $\tilde{GL}$; en fixant la donn\'ee endoscopique elliptique $\underline{G}$, on note $\pi^{GL}_{G,st}$ la projection de cette $\theta$-trace sur $I_{cusp,st}^{G}$  (cf. \ref{pseudocoefficient} (1)).
\begin{thm}Soit $\pi^{GL}$ une repr\'esentation $\theta$-elliptique, irr\'eductible, de $\tilde{GL}$; on suppose que $\pi^{GL}_{G,st}$  n'est ni  nul ni orthogonal \`a toutes les s\'eries discr\`etes de $G$. Alors la $\theta$-trace de $\pi^{GL}$ est un transfert d'un paquet stable de repr\'esentations elliptiques de $G$ ou encore $\pi^{GL}_{H,st}=0$ pour toute donn\'ee endoscopique elliptique $\underline{H}$ diff\'erente de $\underline{G}$.
\end{thm}
On \'ecrit $\pi^{GL}\simeq \times_{(\rho,a)\in {\mathcal E}}St(\rho,a)$ en se souvenant que s'il y a un facteur $F^*$, il op\`ere par $\nu$. Et on pose $\pi^{GL}_{+}:=\times_{(\rho,a)\in {\mathcal E}}St(\rho,a+2) $; c'est une repr\'esentation de $\tilde{GL}(n_{+},E)$, o\`u $n_{+}=n+\sum_{(\rho,a)\in {\mathcal E}}2d_{\rho}$. On d\'ecompose la $\theta$-trace de cette repr\'esentation comme dans \ref{pseudocoefficient} (1). Et on calcule $\circ_{(\rho,a)\in {\mathcal E}}Jac_{\rho\vert\,\vert^{(a+1)/2}}^{GL}$, en ordonnant ${\mathcal E}$ de fa\c{c}on \`a ce que l'on prenne les  $Jac_{\rho\vert\,\vert^{(a+1)/2}}^{GL}$ dans l'ordre croissant sur $a$; \`a chaque \'etape le r\'esultat est une repr\'esentation $\theta$-elliptique irr\'eductible qu'il est facile de calculer. Sa d\'ecomposition dans \ref{pseudocoefficient} (1) se calcule en prenant aussi des modules de Jacquet mais il faut regrouper les donn\'ees endoscopiques qui contiennent le m\^eme sous-groupe de Levi; si $\underline{G}$ est simple, il faut consid\'erer $\underline{G_{n_{+}}}$ mais aussi les donn\'ees compos\'ees de la forme $\underline{G_{1}\times G_{2}}$ tel que $G_{1}= G_{n+2T}$ et $G_{2}$ contient $GL(1/2(n_{+}-n)-T,E)$ comme sous-groupe de Levi maximal.

On v\'erifie que la contribution de ces groupes, si $T\neq (n_{+}-n)/2$,  est nulle: en effet s'il n'en est pas ainsi, on fixe $G_{1}$ et $G_{2}$ pour cette valeur de $T$. Et il existe un paquet stable, $\tau$, de repr\'esentations elliptiques de $G_{2}$ et un sous-ensemble ${\mathcal E}'$ de ${\mathcal E}$ tel que $\circ_{(\rho',a')\in {\mathcal E}'}Jac_{\rho'\vert\,\vert^{(a'+1)/2}}\tau\neq 0$ et:
$$
(n_{+}-n)/2-T=\sum_{(\rho',a')\in {\mathcal E}'}d_{\rho'}.$$ Ceci est impossible car $\tau$ se transf\`ere alors vers une repr\'esentation $\theta$-elliptique, $\pi^{'GL}$, de $\tilde{GL}_{(n_{+}-n)-2T}$ qui doit aussi v\'erifier:
$$
\circ_{(\rho',a')\in {\mathcal E}'}Jac^{GL}_{\rho'\vert\,\vert^{(a'+1)/2}}\pi^{'GL}\neq 0.$$ On sait que $\pi^{'GL}$ est une induite de repr\'esentation de Steinberg. On a donc certainement $$\sum_{(\rho',a')\in {\mathcal E}'}d_{\rho'}(a'+2)\leq (n_{+}-n)-2T$$et ceci n\'ecessite $a'=0$ pour tout $(\rho',a')\in {\mathcal E}'$ ce qui est exclu. D'o\`u l'assertion cherch\'ee.

On note $\pi_{G,st}$ la projection de $\theta$-trace de $\pi^{GL}$ sur $I_{cusp}^{\underline{G}_{n},st}$ et $\pi_{G_{+},st}$ l'analogue pour $G_{n_{+}}$; on sait donc que $\pi_{G,st}$ est un module de Jacquet de $\pi_{G,+}$; soit $\pi_{0}$ une s\'erie discr\`ete de $G$ intervenant dans $\pi_{G,st}$ qui existe par hypoth\`ese. Il existe donc au moins une repr\'esentation temp\'er\'ee $\pi_{0,+}$ intervenant dans $\pi_{+,st}$ et telle que
$$
\circ_{(\rho,a)\in {\mathcal E}}Jac_{\rho\vert\,\vert^{(a+1)/2}}\pi_{0,+}$$ contient $\pi_{0}$ comme sous-quotient; plus pr\'ecis\'ement \`a chaque \'etape dans le calcul du module de Jacquet, il existe un sous-quotient qui est une repr\'esentation temp\'er\'ee et en continuant la proc\'edure, l'\'etape finale donne un sous-quotient qui est $\pi_{0}$;  en fait d'apr\`es \ref{pasdeuxfois}, \`a chaque \'etape le module de Jacquet est une repr\'esentation irr\'eductible et temp\'er\'ee,  la derni\`ere \'etape donnant une s\'erie discr\`ete $\pi_{0}$; ceci n\'ecessite qu'\`a chaque \'etape la repr\'esentation temp\'er\'ee soit en fait une s\'erie discr\`ete. Ainsi $\pi_{0,+}$ est une s\'erie discr\`ete.

On met $\pi_{0,+}$ dans un paquet stable comme en \ref{projectionnonnulle}; on transfert ce paquet stable en la $\theta$-trace d'une repr\'esentation virtuelle combinaison lin\'eaire de repr\'esentations $\theta$-elliptiques de $\tilde{GL}(n_{+},E)$, not\'e $\tilde{\pi}^{GL}_{0,+}$ et on recalcule $$\circ_{(\rho,a)\in {\mathcal E}}Jac_{\rho\vert\,\vert^{(a+1)/2}}\tilde{\pi}^{GL}_{0,+};$$ le r\'esultat est un transfert d'une repr\'esentation virtuelle stable combinaison lin\'eaire de s\'eries discr\`etes de $G_{n}$ qui n'est certainement pas nulle car elle contient n\'ecessairement $\pi_{0}$ avec un coefficient non  nul (cf. \ref{pasdeuxfois} (ii)). Mais comme
$$
\sum_{(\rho,a)\in {\mathcal E}} d_{\rho}(a+2)=n_{+},$$
n\'ecessairement $\tilde{\pi}^{GL}_{0,+}$ contient l'induite $(\times_{(\rho,a)\in {\mathcal E}}St(\rho,a+2))\otimes \nu$ \`a un scalaire pr\`es et le module de Jacquet consid\'er\'e est $\pi^{GL}$ (celui de l'\'enonc\'e) \`a un scalaire pr\`es. Ainsi $\pi^{GL}$ est un transfert d'une distribution stable pour $\underline{G}_{n}$, ce qui est le r\'esultat cherch\'e.

\section{Le cas des repr\'esentations cuspidales}

\subsection{Point de r\'eductibilit\'e des induites de cuspidales\label{redcusp}}
Soit $\pi$ une repr\'esentation cuspidale irr\'eductible de $G$ et soit $\rho$ une repr\'esen\-tation cuspidale irr\'eductible et unitaire de $GL(d_{\rho},F)$ ce qui d\'efinit $d_{\rho}$. 
\begin{thm}Soit $x\in {\mathbb R}$ tel que l'induite de la repr\'esentation $\rho\vert\,\vert^{x}\times \pi$, qui est une repr\'esentation de $G_{n+d_{\rho}}$ soit r\'eductible, alors $x\in 1/2 {\mathbb Z}$ et $x\leq (n/d_{\rho}+1)/2$.
\end{thm}
On fixe $\pi,\rho,x$ comme dans l'\'enonc\'e; si $x=0$, le th\'eor\`eme est clair. Si l'induite temp\'er\'ee $\rho\times \pi$ est sans entrelacement sous l'action du groupe de Weyl, d'apr\`es Harish-Chandra, $x$ ne peut exister. On suppose donc que cette induite a un entrelacement et cela permet de supposer que $x>0$. On note alors $\pi_{+}$ l'unique sous-module irr\'eductible de l'induite $\rho\vert\,\vert^{x}\times \pi$; c'est une s\'erie discr\`ete. Elle appartient \`a un paquet stable d'apr\`es \ref{projectionnonnulle}, not\'e $\pi_{+,st}$ et il existe donc une repr\'esentation virtuelle not\'ee $\pi^{GL}_{+}$ dont la trace tordue par $\theta$ est un transfert de $\pi_{+,st}$. 

On sait que $Jac_{\rho\vert\,\vert^x}^{GL}\pi^{GL}_{+}$ est un transfert de $Jac_{\rho\vert\,\vert^x}\pi_{+,st}$ du paquet stable pr\'ec\'edent (cf la preuve de \ref{projectionendoscopique} ); ce dernier module de Jacquet contient de fa\c{c}on non trivial la trace de $\pi$ car il n'y a qu'une repr\'esentation irr\'eductible de $G_{n}$ qui a dans son module de Jacquet cuspidal $\rho\vert\,\vert^{x}\otimes \pi$ et c'est $\pi_{+}$. Ainsi la repr\'esentation virtuelle obtenue par module de Jacquet n'est pas nulle donc son transfert n'est pas nulle. Cela force aussi l'in\'egalit\'e de l'\'enonc\'e car on a n\'ecessairement $n+d_{\rho}\geq d_{\rho}(2x+1)$.

On sait que $\pi^{GL}_{+}$ est une combinaison lin\'eaire d'induites de repr\'esentations de Steinberg unitaire; son module de Jacquet $Jac_{\rho\vert\,\vert^{x}}^{GL}(\pi^{GL}_{+})$ est non nul et cela  force $x$ \`a \^etre un demi-entier. Cela prouve le th\'eor\`eme.
\subsection{Blocs de Jordan des repr\'esentations cuspidales\label{cascuspidal}}
Soit $\pi$ une repr\'esentation cuspidale irr\'eductible de $G_{n}$.
On d\'efinit ici $Jord(\pi)$ comme \'etant l'ensemble des couples $(\rho,a)$ o\`u $\rho$ est une repr\'esentation cuspidale de $GL(d_{\rho},F)$ et  o\`u $a$ est un entier, $Jord(\pi)$ satisfaisant \`a

 si $(\rho,a)\in Jord(\pi)$ avec $a>2$ alors $(\rho,a-2)\in Jord(\pi)$;
 
  $(\rho,a)\in Jord(\pi)$ (avec $a\geq 1$) et pour tout $a'>a$, $(\rho,a')\notin Jord(\pi)$ si et seulement si l'induite $\rho\vert\,\vert^{(a+1)/2}\times \pi$ de $G_{n^*+d_{\rho}}$ est r\'eductible.

On d\'emontrera en \ref{unicitedefinition} que cette d\'efinition correspond aux d\'efinitions de \cite{europe} et \cite{jams}.
\begin{thm} L'ensemble $Jord(\pi)$ est fini et $\sum_{(\rho,a)\in Jord(\pi)}ad_{\rho}=n$. On pose $$\pi^{GL}:=\times_{(\rho,a)\in Jord(\pi)}St(\rho,a),$$ cette repr\'esen\-tation est $\theta$-invariante et sa trace tordue est un transfert d'une combinaison lin\'eaire finie de traces de s\'eries discr\`etes de $G$, l'une d'entre elles \'etant $\pi$ (avec un coefficient non nul).
 
De plus toute combinaison lin\'eaire stable de repr\'esentations elliptiques de $G$ contenant $\pi$ avec un coefficient non nul et se transf\'erant sur la $\theta$-trace d'une repr\'esentation irr\'eductible et $\theta$-elliptique de $\tilde{GL}$ se transf\`ere n\'ecessairement sur la repr\'esentation $\pi^{GL}$ qui vient d'\^etre d\'ecrite (\`a un scalaire pr\`es).
 \end{thm}
 Avant de faire la d\'emonstration remarquons que ce th\'eor\`eme n'est pas tout \`a fait ce que l'on veut et c'est ce qui explique sa formulation un peu compliqu\'ee; on veut d\'emontrer et on d\'emontrera que la combinaison lin\'eaire stable dont il est question au d\'ebut est celle qui est associ\'ee \`a $\pi$ par \ref{projectionnonnulle}.

Soit ${\mathcal E}$ un sous-ensemble fini de $Jord(\pi)$ tel que si $(\rho,a)\in {\mathcal E}$ et $(\rho,a+2)\in Jord(\pi)$ alors $(\rho,a+2)\in {\mathcal E}$. On v\'erifie la propri\'et\'e suivante:

il existe une  repr\'esentation irr\'eductible de $G_{(n+\sum_{(\rho,a)\in {\mathcal E}}2d_{\rho},F)}$, not\'ee $\pi_{\mathcal E}$ telle que la repr\'esentation virtuelle
$$\circ_{(\rho,a)\in {\mathcal E}}Jac_{\rho\vert\,\vert^{(a+1)/2}}\pi_{\mathcal E}$$ contient $\pi$ (avec un coefficient non nul), o\`u ${\mathcal E}$ est muni d'un ordre tel que l'on prend $Jac_{\rho\vert\,\vert^{(a+1)/2}}$ avant $Jac_{\rho\vert\,\vert^{(a'+1)/2}}$ si $a>a'$. De plus $\pi_{\mathcal E}$ est une s\'erie discr\`ete.

Pour cela, on consid\`ere la repr\'esentation induite
$$
\times_{(\rho,a)\in {\mathcal E}}\rho\vert\,\vert^{(a+1)/2}\times \pi.
$$
On v\'erifie que cette induite a un unique sous-module irr\'eductible en utilisant la r\'eciprocit\'e de Frobenius. Ce sous-module irr\'eductible est not\'e  $\pi_{{\mathcal E}}$. Il est facile de calculer le module de Jacquet cuspidal de cette sous-repr\'esentation:  c'est la somme directe des termes $\otimes_{i\in [1,\vert{\mathcal E}\vert]}\rho_{i}\vert\,\vert^{(a_{i}+1)}\otimes \pi,$ o\`u  $\{(\rho_{i},a_{i}) ; i \in [1,\vert{\mathcal E}\vert]\}$ est exactement l'ensemble ${\mathcal E}$ ordonn\'e de telle sorte que si pour $i<j$, $\rho_{i}=\rho_{j}$ alors $a_{i}<a_{j}$. Il est alors clair que $\pi_{{\mathcal E}}$ est une s\'erie discr\`ete.
 
On pose $n_{+}=n+\sum_{(\rho,a)\in {\mathcal E}}2d_{\rho}$ et on consid\`ere une repr\'esentation virtuelle, $\pi^{GL}_{\mathcal E}$ de $\tilde{GL}(n_{+},E)$ $\theta$ invariante qui est un transfert du paquet stable contenant $\pi_{\mathcal E}$ tel que d\'efini en \ref{projectionnonnulle}. On calcule le module de Jacquet $ \circ_{(\rho,a)\in {\mathcal E}}Jac^{GL}_{\rho\vert\,\vert^{(a+1)/2}}$ comme on l'a fait ci-dessus et on conclut que le r\'esultat est non nul car c'est un transfert d'une repr\'esentation virtuelle contenant $\pi$ avec un coefficient non nul. Ainsi $\pi^{GL}_{\mathcal E}$, qui rappelons-le est une repr\'esentation virtuelle, contient au moins une repr\'esentation induite de repr\'esentations de Steinberg v\'erifiant la non nullit\'e ci-dessus; il est facile de calculer les modules de Jacquet des repr\'esentations de Steinberg. En prenant le premier $Jac^{GL}_{\rho_{0}\vert\,\vert^{(a_{0}+1)/2}}$ on voit  que cette repr\'esenta\-tion est n\'ecessairement de la forme une induite de repr\'esentation $St(\rho_{0},a_{0}+2)$  avec \'eventuellement une autre repr\'esentation temp\'er\'ee, $\tau$. Le r\'esultat est alors l'induite de $St(\rho,a)$ avec $\tau$ et on doit encore avoir la non nullit\'e $$\circ_{(\rho,a)\in {\mathcal E}- \{(\rho_{0},a_{0})}Jac^{GL}_{\rho\vert\,\vert^{(a+1)/2}}\neq 0;$$ l'hypoth\`ese sur l'ordre assure que dans l'ensemble restant si $\rho\simeq \rho_{0}$ alors $a>a_{0}$. La non nullit\'e porte alors sur le mo\-dule de Jacquet de $\tau$ et on en proc\'edant ainsi on voit que la repr\'esentation fix\'ee est une induite de $\times_{(\rho,a)\in {\mathcal E}}St(\rho,a)$ avec \'eventuellement encore une repr\'esentation temp\'er\'ee.  Mais on a donc n\'ecessaire\-ment
 $
 \sum_{(\rho,a)\in {\mathcal E}}d_{\rho}(a+2) \leq n+\sum_{(\rho,a)\in {\mathcal E}}2 d_{\rho}.
 $
D'o\`u $\sum_{(\rho,a)\in {\mathcal E}}ad_{\rho}\leq n$; comme $n$ est ind\'ependant de ${\mathcal E}$, cela force la finitude de $Jord(\pi)$ et l'in\'egalit\'e
 $$
 \sum_{(\rho,a)\in Jord(\pi)}ad_{\rho}\leq n.$$
 On pose $n_{-}:= n-\sum_{(\rho,a)\in Jord(\pi)}ad_{\rho}$ et on vient de montrer qu'il existe une repr\'esentation virutelle $\tau_{0}$ de $\tilde{GL}(n_{-},E)$ tel que 
 $\times_{(\rho,a)\in Jord(\pi)}St(\rho,a)\times \tau_{0}$ soit $\theta$-invariante et soit un transfert du paquet stable associ\'e \`a $\pi$ dans \ref{projectionnonnulle}.

 Pour montrer que $n_{-}=0$, on proc\`ede de fa\c{c}on inverse: on fixe $\pi^{GL}$ une repr\'esentation irr\'eductible et $\theta$-elliptique de $\tilde{GL}(n)$ dont la projection de la $\theta$-trace sur $I_{cusp}^{\underline{G},st}$ est non orthogonale \`a la trace de $\pi$ (cf. \ref{pseudocoefficient} (1)). On d\'efinit $Jord({\pi}^{GL})$ comme \'etant l'ensemble des couples $(\rho,a)$  tel que 
 $$
 \pi^{GL}=\times_{(\rho,a)\in Jord(\pi^{GL})}St(\rho,a).
 $$
 On pose $n_{+}:= n+\sum_{(\rho,a)\in Jord(\pi^{GL})2d_{\rho}}$ et on  consid\`ere $\pi^{GL}_{+}$ la repr\'esentation irr\'eductible du groupe $\tilde{GL}(n_{+})$ \'egale \`a $\times_{(\rho,a)\in Jord(\pi^{GL})}St(\rho,a+2)$. On note $\tau_{+}$ la repr\'esentation virtuelle de $G_{n_{+}}$ obtenu par la projection de la  $\theta$-trace de cette repr\'esentation sur $I_{cusp}^{\underline{G}_{n_{+}}}$. 
 On remarque qu'en ordonnant $Jord({\pi}^{GL}_{+})$ de sorte que on prenne d'abord $Jac_{\rho\vert\,\vert^{(a+1)/2}}^{GL}$ avant de prendre $Jac_{\rho\vert\,\vert^{(a'+1)/2}}^{GL}
 $ pour le m\^eme $\rho$ et la m\^eme parit\'e de $a,a'$ si  $a'>a$, on a
 $$
 \circ_{(\rho,a)\in Jord(\pi^{GL})} Jac^{GL}_{\rho\vert\,\vert{(a+1)/2}}\pi^{GL}_{+}=\pi^{GL}.
 $$
 C'est le m\^eme argument que celui donn\'e ci-dessus.
  Puisque le transfert commute au module de jacquet, on a s\^urement que 
$\circ_{(\rho,a)\in Jord(\pi^{GL})} Jac_{\rho\vert\,\vert^{(a+1)/2}}\tau_{+}$ contient $\pi$. Ainsi il existe $\pi_{+}$ une s\'erie discr\`ete intervenant dans $\tau_{+}$ et tel qu'\`a un scalaire pr\`es, $
 \circ_{(\rho,a)\in Jord(\pi^{GL})} Jac_{\rho\vert\,\vert{(a+1)/2}}\pi_{+}=\pi;
 $ on a m\^eme un r\'esultat un peu plus pr\'ecis que l'on n'exprime que dans le cas qui nous int\'eresse: on fixe $\rho_{0}$ une repr\'esentation cuspidale unitaire et irr\'eductible, une parit\'e et $a_{0}$ maximal avec cette parit\'e fix\'ee  tels que $Jord(\pi^{GL})$ contienne $(\rho_{0},a_{0})$. On note $Jord(\pi^{GL})_{-}$ l'ensemble $Jord(\pi^{GL})$ priv\'e de $(\rho_{0},a_{0})$ et on a:
 $$
  \circ_{(\rho,a)\in Jord(\pi^{GL})} Jac_{\rho\vert\,\vert^{(a+1)/2}}\tau_{+}= Jac^{GL}_{\rho_{0}\vert\,\vert^{(a_{0}+1)/2}}$$
  $$\circ_{(\rho,a)\in Jord(\pi^{GL})_{-}} Jac_{\rho\vert\,\vert^{(a+1)/2}}\tau_{+}.$$
Or $\circ_{(\rho,a)\in Jord(\pi^{GL})_{-}} Jac_{\rho\vert\,\vert^{(a+1)/2}}\tau_{+}$ n'est autre que $$St(\rho_{0},a_{0}+2)\times_{(\rho,a)\in Jord(\pi^{GL})_{-}}St(\rho,a).$$ Comme ci-dessus, cela montre qu'il existe une s\'erie discr\`ete $\pi_{0}$ de $G_{n+2d_{\rho_{0}}}$ tel que $Jac_{\rho_{0}\vert\,\vert^{(a_{0}+1)/2}}\pi_{0}$ contient $\pi$. Ainsi $\pi_{0}$ est un sous-quotient de l'induite $\rho\vert\,\vert^{(a_{0}+1)/2}\times \pi$ et cette induite est r\'eductible. Cela entra\^{\i}ne par d\'efinition que $(\rho_{0},a_{0})\in Jord(\pi)$ et par les propri\'et\'es de $Jord(\pi)$ pour tout $a'$ de m\^eme parit\'e que $a_{0}$ si $(\rho_{0},a')\in Jord(\pi^{GL})$ alors $(\rho_{0},a')\in Jord(\pi)$. Ainsi on a 
$$
n=\sum_{(\rho,a)\in Jord(\pi^{GL})}ad_{\rho}\leq \sum_{(\rho,a)\in Jord(\pi)}ad_{\rho}=n-n_{-}.
$$
Cela force $n_{-}=0$ et  $Jord(\pi^{GL})=Jord(\pi)$. En revenant au d\'ebut de la preuve, on a en plus montr\'e que $\times_{(\rho,a)\in Jord(\pi)}St(\rho,a)$ est un transfert d'un \'el\'ement de $I_{cusp}^{\underline{G},st}$ non orthogonal \`a $\pi$. Mais, puisque  $Jord(\pi^{GL})=Jord(\pi)$, 
 $\times_{(\rho,a)\in Jord(\pi)}St(\rho,a)=\pi^{GL}$ et il n'y a donc qu'une repr\'esentation $\theta$-elliptique et irr\'eductible de $\pi^{GL}$ dont la projection de la $\theta$-trace sur $I_{cusp}^{\underline{G},st}$ n'est pas orthogonale \`a $\pi$.
\section{Support cuspidal \'etendu\label{supportcuspidaletendu}}
En \ref{cascuspidal}, pour $\pi$ une repr\'esentation cuspidale, on a d\'efini $Jord(\pi)$ et $\pi^{GL}$; on appelle support cuspidal \'etendu de $\pi$ le support cuspidal de $\pi^{GL}$ et on voit ce support cuspidal comme un ensemble de couples $(\rho,x)$, o\`u $\rho$ est une repr\'esentation cuspidale unitaire d'un groupe $GL(d_{\rho},E)$ et $x$ est un nombre r\'eel en fait un demi-entier relatif. Cet ensemble est bien d\'efini \`a permutation pr\`es.

On g\'en\'eralise cette d\'efinition \`a toute repr\'esentation irr\'eductible:

soit $\pi$ une repr\'esentation irr\'eductible de $\underline{G}$; on \'ecrit $\pi$ comme sous-quotient d'une repr\'esentation induite \`a partir d'une repr\'esentation cuspidale d'un sous-groupe de Levi. Cela d\'efinit le support cuspidal usuel de $\pi$ que l'on \'ecrit comme une collection de couples $(\rho,x)$ et d'une repr\'esentation cuspidale $\pi_{cusp}$ d'un groupe de m\^eme type que $\underline{G}$ mais en g\'en\'eral de rang plus petit, o\`u les $(\rho,x)$ sont des couples form\'es d'une repr\'esentation cuspidale unitaire $\rho$ et d'un nombre r\'eel $x$. Alors que $\pi_{cusp}$ est uniquement d\'etermin\'e, les couples $(\rho,x)$ sont d\'efinis \`a permutation pr\`es et \`a changement de $(\rho,x)$ en $(\tilde{\rho},-x)$ o\`u, ici, $\tilde{\rho}:=\rho^*\otimes \nu$ si $E=F$ et $\rho{\check\empty}$ si $E\neq F$ o\`u $\nu$ est le caract\`ere de \ref{groupes}; c'est un calcul sur le groupe de Weyl. On d\'efinit le support cuspidal \'etendu de $\pi$ comme \'etant l'union du support cuspidal \'etendu de $\pi_{cusp}$ avec l'ensemble des couples $(\rho,x), (\tilde{\rho},-x)$ o\`u $(\rho,x)$ parcourt l'ensemble ci-dessus. Le support cuspidal \'etendu et donc bien d\'efini \`a permutation pr\`es. 

\begin{rmq} Le support cuspidal \'etendu de $\pi$ est n\'ecessairement le support cuspidal d'au moins une repr\'esentation irr\'eductible de $GL(n,E)$ et d'au plus une repr\'esentation temp\'er\'ee de ce groupe.
\end{rmq}
La premi\`ere assertion vient uniquement d'un calcul de dimension: on note $Supp(\pi)$ le support cuspidal de $\pi$ sauf $\pi_{cusp}$ (on fait un choix qui n'aura pas d'importance), on note $Supp_{et}(\pi)$ et $Supp_{et}(\pi_{cusp})$ le support cuspidal \'etendu de $\pi$ et $\pi_{cusp}$ et on a:
$$
\sum_{(\rho,x)\in Supp_{et}(\pi)}d_{\rho}=2\sum_{(\rho,x)\in Supp(\pi)}d_{\rho}+\sum_{(\rho,x)\in Supp_{et}(\pi_{cusp})}d_{\rho}=n.
$$
La deuxi\`eme assertion de la remarque vient du fait que le support cuspidal d'une repr\'esentation de $GL(n,E)$ est une union de segments centr\'es en $0$ et que la repr\'esentation temp\'er\'ee est uniquement d\'etermin\'ee par cette union de segments; mais il n'y a au plus qu'une fa\c{c}on d'\'ecrire un ensemble de repr\'esentations cuspidales comme union de segments centr\'es en $0$, d'o\`u l'assertion.

\

Le but de la suite du travail est de montrer que si $\pi$ est une repr\'esentation temp\'er\'ee, alors le support cuspidal \'etendu de $\pi$ est le support cuspidal d'une repr\'esentation temp\'er\'ee et $\theta$-invariante de $GL(n,E)$ et que cette repr\'esentation temp\'er\'ee a sa $\theta$-trace qui est un transfert d'un paquet stable (en fait ''du'' si ceci est bien d\'efini) de $G$ contenant $\pi$.

\subsection{Le cas des s\'eries discr\`etes strictement positives\label{seriesdiscretespositives}}
 Une g\'en\'eralisation des repr\'esentations cuspidales, pour ce qui est fait ici, est la notion de s\'eries discr\`etes strictement positives, comme cela avait d\'ej\`a \'et\'e le cas en \cite{europe}. Par d\'efinition une telle s\'erie discr\`ete est une repr\'esentation irr\'eductible dont les modules de Jacquet cuspidaux sont de la forme $\otimes_{(\rho,x)}\rho\vert\,\vert^x\times \pi_{cusp}$, o\`u $(\rho,x)$ d\'ecrit un ensemble de couples form\'es d'une repr\'esentation cuspidale unitaire et d'un nombre r\'eel $x$ strictement positif et o\`u $\pi_{cusp}$ est une repr\'esentation cuspidale d'un groupe de la forme $G_{n'}$ avec $n'\leq n$. Dans le cas o\`u $G=SO(2n,F)$ o\`u $GSpin_{2n}(F)$, il faut consid\'erer simultan\'ement $\pi$ et 
 $\tilde{\pi}$ l'image de $\pi$ par l'automorphisme ext\'erieur et demander la m\^eme propri\'et\'e pour ces deux repr\'esentations, ou encore  ceci revient \`a travailler avec le groupe non connexe $O(2n,F)$ et son analogue pour $GSpin_{2n}(F)$.  Ici on ne classifiera pas ces repr\'esentations (cf. \cite{europe} qui se g\'en\'eralise cf par exemple \cite{kim407} et Jang Bogume \cite{bogume}), les propri\'et\'es d\'ecrites ci-dessous nous suffisant.
 
Soit $(\rho,x)$ un couple form\'e d'une repr\'esentation cuspidale unitaire  $\rho$  et d'un nombre r\'eel strictement positif $x$. Et soit $\pi$ une s\'erie discr\`ete fortement positive.

\begin{lem} (i) On suppose que $Jac_{\rho\vert\,\vert^x}\pi\neq 0$; alors $Jac_{\rho\vert\,\vert^x}\pi$ est une s\'erie discr\`ete strictement positive.

(ii)Supposons que $x\geq 1$; les sous-quotients irr\'eductibles de l'induite $\rho\vert\,\vert^x\times \pi$ sont  des s\'eries discr\`etes fortement positives sauf l'unique quotient irr\'eductible qui peut \'eventuellement \^etre toute l'induite.
\end{lem}
Le (i) est compl\`etement \'evident. 
Montrons  (ii); soit $\pi'$ un sous-quotient irr\'educ\-tible de l'induite de l'\'enonc\'e. Les modules de Jacquet cuspidaux de l'induite sont obtenus en prenant ceux de $\pi$ et en glissant soit $\rho\vert\,\vert^x$ soit $\theta(\rho)\vert\,\vert^{-x}$ \`a n'importe quelle place sauf si $G=SO(2n,F)$ ou $GSpin_{2n}(F)$, o\`u il faut consid\'erer simultan\'ement $\pi$ et $\tilde{\pi}$ (l'image de $\pi$ par l'automorphisme ext\'erieur). Supposons que $\pi'$ ait un module de Jacquet contenant un terme avec $\theta(\rho)\vert\,\vert^{-x}$. Par r\'eciprocit\'e de Frobenius on trouve une inclusion:
$$
\pi'\hookrightarrow \times_{(\rho',x')}\rho'\vert\,\vert^{x'}\times \theta(\rho)\vert\,\vert^{-x}\times \tau,
$$
o\`u $\tau$ est une repr\'esentation irr\'eductible convenable; n\'ecessairement ici $x'\geq 1/2$ d'o\`u $x+x'>1$ par l'hypoth\`ese sur $x'$ et on peut mettre $\theta(\rho)\vert\,\vert^{-x}$ en premi\`ere position. Et $\pi'$ est alors le quotient de Langlands de l'induite.
\begin{cor} Il existe une unique repr\'esentation $\theta$-elliptique  de $\tilde{GL}$ qui est irr\'e\-ductible et  dont la $\theta$-trace est un transfert d'une combinaison lin\'eaire stable  de repr\'esentations elliptiques de $G$ contenant $\pi$. De plus le support cuspidal de cette repr\'esentation de ${GL}(n,E)$ est exactement le support cuspidal \'etendu de $\pi$.
\end{cor}
L'existence de cette repr\'esentation r\'esulte de \ref{pseudocoefficient} (1) et de \ref{projectionendoscopique}.
On d\'emontre l'unicit\'e et sa description par r\'ecurrence sur $n$ en initialisant la r\'ecurrence avec les repr\'esentations cuspidales si $n>0$; et dans le cas o\`u $n=0$, le r\'esultat est vide sauf pour $SO(1,F)$ et $GSpin(1,F)$ o\`u il est trivial, toutefois, il faut dans ce cas accept\'e une petite extension de la formulation du th\'eor\`eme pour que $\tilde{GL}$ existe aussi avec $n=0$.

Soit $\pi$ une s\'erie discr\`ete strictement positive et soit $\rho,x$ comme dans le lemme (i) ci-dessus. Soit $\pi^{GL}$ une repr\'esentation $\theta$-elliptique irr\'eductible de $\tilde{GL}$ qui est un transfert d'une combinaison lin\'eaire stable  de repr\'esentations elliptiques de $G$ contenant $\pi$ (cf. \ref{pseudocoefficient} (1) et \ref{projectionendoscopique}).  Il suffit \'evidemment de d\'emontrer que le support cuspidal de $\pi^{GL}$ co\"{\i}ncide avec le support cuspidal \'etendu de $\pi$ et on peut supposer que $\pi$ n'est pas cuspidal. On fixe $(\rho,x)$ tel que $Jac_{\rho\vert\,\vert^x}\pi\neq 0$ et on applique l'hypoth\`ese de r\'ecurrence \`a $\pi_{-}:=Jac_{\rho\vert\,\vert^x}\pi$. On note $\tilde{\pi}_{st}$ la combinaison lin\'eaire stable  de repr\'esentations elliptiques de $G$ contenant $\pi$  se transf\'erant en la $\theta$-trace de $\pi^{GL}$. La repr\'esentation $Jac_{\rho\vert\,\vert^{x}}^{GL}(\pi^{GL})$ est un transfert de $Jac_{\rho\vert\,\vert^{x}}\tilde{\pi}_{st}$; comme cette derni\`ere repr\'esentation contient une s\'erie discr\`ete, $\pi_{-}$, la repr\'esentation $$Jac_{\rho\vert\,\vert^{x}}\tilde{\pi}_{st}$$ ne peut pas \^etre une induite $\theta$-stable \`a partir d'un espace de Levi de $\tilde{GL}$. Ainsi si l'on \'ecrit $\pi^{GL}$ comme l'induite $\times_{(\rho',a')\in {\mathcal E}}St(\rho',a')$, n\'ecessairement $(\rho,2x+1)\in {\mathcal E}$ et $(\rho,2x-3)\notin {\mathcal E}$; deplus, $St(\rho,a-2)\times_{(\rho',a')\in {\mathcal E}; (\rho',a')\neq (\rho,a)}St(\rho',a')$ est un transfert n\'ecessairement d'une combinaison lin\'eaire de repr\'esentations elliptiques de $G_{n-2d_{\rho}}$ ($d_{\rho}$ est d\'efini par le fait que $\rho$ est une repr\'esentation de $GL(d_{\rho},E)$) et qui co\"{\i}ncide avec $Jac_{\rho\vert\,\vert^x}\tilde{\pi}_{st}$. Cette combinaison lin\'eaire contient donc $\pi_{-}$. Par l'hypoth\`ese de r\'ecurrence, le support cuspidal \'etendu de $\pi_{-}$ est exactement le support cuspidal de $Jac_{\rho\vert\,\vert^x}^{GL}(\pi^{GL})$. Le support cuspidal de $\pi$ s'obtient en rajoutant $\rho\vert\,\vert^x$ et $\rho\vert\,\vert^{-x}$ au support cuspidal \'etendu de $\pi_{-}$; c'est la m\^eme op\'eration qui fait passer du support cuspidal de  de $Jac_{\rho\vert\,\vert^{x}}^{GL}\pi^{GL}$ \`a $\pi^{GL}$ d'o\`u le r\'esultat.
\subsection{Bloc de Jordan des s\'eries discr\`etes strictement positives\label{jordpositives}}
Soit $\pi$ une s\'erie discr\`ete strictement positive; on lui a associ\'e dans le corollaire de  \ref{seriesdiscretespositives} une unique repr\'esentation temp\'er\'ee $\pi^{GL}$ de $\tilde{GL}(n,E)$. Ici on prend comme d\'efinition de $Jord(\pi):=\{(\rho,a)\}$ de  telle sorte que $$\pi^{GL}\simeq \times_{(\rho,a)\in Jord(\pi)}St(\rho,a);$$ on v\'erifiera dans la remarque de  \ref{cardinal} que cette d\'efinition  est bien analogue \`a celle donn\'ee en \cite{pourshahidi}. Mais ici on d\'emontre le r\'esultat technique dont nous aurons besoin: soit $\ell \in {\mathbb N}$ et pour tout $i\in [1,\ell]$ la donn\'ee d'un couple $\rho_{i},x_{i}\geq 1$ et d'une s\'erie discr\`ete strictement positive $\pi_{i}$ v\'erifiant inductivement que $\pi_{i}$ est un sous-quotient de l'induite $\rho_{i}\vert\,\vert^{x_{i}}\times \pi_{i-1}$, o\`u l'on a pos\'e $\pi_{0}=\pi$ et $x_{i}>x_{i-1}$ avec $x_{0}=1/2$.
\begin{lem} Pour une telle suite, l'ensemble $\{(\rho_{i},2x_{i}-1); i\in [1,\ell]\}$ est un sous-ensemble de $Jord(\pi)$.
\end{lem}
On fixe $i\in [1,\ell]$; on admet par r\'ecurrence que $Jord(\pi_{i-1})$ s'obtient \`a partir de $Jord(\pi)$ en rempla\c{c}ant les blocs $(\rho_{j},2x_{j}-1); j\in [1,i[$ par les blocs $(\rho_{j},2x_{j}+1)$, cette assertion est vraie pour $i=0$ et on la d\'emontre pour le $i$ fix\'e; on sait que le support cuspidal \'etendu de $\pi_{i}$ est celui de $\pi_{i-1}$ auquel on ajoute $(\rho_{i},x_{i}), (\tilde{\rho}_{i},-x_{i})$ (o\`u $\tilde{\rho}\simeq \rho^*\otimes \nu$ ou $\rho^{\check{\empty}}$ (cf. \ref{groupes})). Mais d'autre part, comme $\pi_{i}$ est par hypoth\`ese une s\'erie discr\`ete strictement positive, son support cuspidal \'etendu est une union de segment; comme $x_{i}>1/2$ la seule possibilit\'e est que $\rho_{i}\simeq \tilde{\rho_{i}}$ et que l'on rajoute $(\rho_{i},x_{i})$ \`a l'extr\'emit\'e d'un des segments du support cuspidal \'etendu de $\pi_{i-1}$; il faut donc que  $(\rho_{i},2x_{i}-1)\in Jord(\pi_{i-1})$. On utilise alors l'hypoth\`ese que $x_{i}<x_{j}$ pour tout $j\in [1,i[$ ce qui force que $(\rho_{i},2x_{i}-1)\in Jord(\pi)$ et c'est alors ce couple qui devient $(\rho_{i},2x_{i}+1)$ dans $Jord(\pi_{i})$. D'o\`u le lemme.

\subsection{Un lemme de structure des s\'eries discr\`etes et des repr\'esentations temp\'er\'ees\label{lemmedestructure}}
Soit $\rho$ une repr\'esentation cuspidale de $GL(d_{\rho},E)$. Soit $d,f$ des nombres r\'eels tel que $d+f\in {\mathbb N}_{\geq 0}$. On note $
\langle d,-f\rangle _{\rho}$ l'unique sous-repr\'esentation de $GL((d+f+1)d_{\rho},F)$ incluse dans l'induite $\times_{i\in [d,-f]}\rho\vert\,\vert^{i}$; c'est une s\'erie discr\`ete tordu par un caract\`ere non unitaire.
\begin{lem} Soit $\pi$ une repr\'esentation irr\'eductible de $G$. Il existe un ensemble, ${\mathcal T}$ de triplets $(\rho,d,f)$ ordonn\'e  et une s\'erie discr\`ete strictement positive, $\pi_{>0}$ tels que l'on ait une inclusion
$$
\pi \hookrightarrow \times _{(\rho,d,f)\in {\mathcal T}} \langle d,-f\rangle _{\rho} \times \pi_{>0}
$$ et tels que $ f\geq 0$ pour tout $(\rho,d,f)\in {\mathcal T}$ et $f\geq f'$ si le triplet $(\rho,d,f)$ pr\'ec\`ede le triplet $(\rho',d',f')$ dans l'ordre de ${\mathcal T}$.
\end{lem}
Pr\'ecisons la situation des groupes $SO(2n,F)$ et $GSpin_{2n}(F)$; on ne voit les repr\'esentations qu'\`a conjugaison pr\`es par l'automorphisme ext\'erieur; donc dans le lemme, il faut \'eventuellement remplacer $\pi$ par son conjugu\'e par cet automorphisme.

On fixe une inclusion de $\pi$ dans une repr\'esentation induite \`a partir d'une repr\'esentation cuspidale: ci-dessous $\ell\in {\mathbb N}$, pour $i\in [1,\ell]$, les $\rho_{i}$ sont des repr\'esentations cuspidales unitaires, les $x_{i}$ sont des nombres r\'eels et $\pi_{cusp}$ est une repr\'esentation cuspidale d'un groupe $G_{n_{cusp}}$, pour $n_{cusp}$ convenable:
$$
\pi\hookrightarrow \times_{i\in [1,\ell]}\rho_{i}\vert\,\vert^{x_{i}}\times \pi_{cusp}; \eqno(1)
$$ Le terme de droite de (1) n'est \'evidemment pas uniquement d\'etermin\'e et on fait un choix de telle sorte que le nombre de $min_{i\in [1,\ell]}x_{i}$  soit minimal pour tous les choix possibles. Ainsi si dans (1) aucun des $x_{i}$ n'est inf\'erieur ou \'egal \`a $0$, $\pi$ est une s\'erie discr\`ete strictement positive (par r\'eciprocit\'e de Frobenius aucun terme du module de Jacquet cuspidal ne contient de termes ''n\'egatifs'').

On suppose donc qu'il existe $x_{i}\leq 0$ dans le terme de droite de (1) ayant la propri\'et\'e de minimalit\'e pr\'ec\'edente et on note $\ell_{1}$ le plus petit indice tel que $x_{\ell_{1}}$ soit minimal dans (1). On pose $f_{1}:=-x_{\ell_{1}}$; on peut modifier le terme de droite de (1) uniquement en permutant \'eventuellement les termes $x_{i}$ pour $i\leq \ell_{1}$ de sorte que $\rho_{\ell_{1}}\vert\,\vert^{x_{\ell_{1}}}$ soit le plus \`a gauche possible. On est alors s\^ur que pour tout $i\leq \ell_{1}$, $\rho_{i}\simeq \rho_{\ell_{1}}$ et il existe un sous-quotient $\tau_{1}$ de l'induite de $GL(d_{\rho_{\ell_{1}}}\ell_{1},E)$ $\times_{i\in [1,\ell_{1}]}\rho_{\ell_{1}}\vert\,\vert^{x_{1}}$ tel que l'inclusion (1) se factorise en 
$$
\pi\hookrightarrow \tau_{1}\times_{i>\ell_{1}}\rho_{i}\vert\,\vert^{x_{i}}\times \pi_{cusp}. \eqno(2)
$$
La minimalit\'e de $\ell_{1}$ assure alors que $\tau_{1}= \langle x_{1},-f_{1}\rangle_{\rho_{\ell_{1}}}$. On peut encore factoriser l'inclusion (2) en
$$
\pi\hookrightarrow \tau_{1}\times \pi',
$$
o\`u $\pi'$ est une repr\'esentation convenable. On applique le lemme \`a $\pi'$ et on obtient le lemme pour $\pi$ en remarquant que l'in\'egalit\'e du lemme est bien satisfaite par minimalit\'e de $x_{\ell_{1}}$.
\begin{rmq} On suppose que $\pi$ est une repr\'esentation temp\'er\'ee (resp. une s\'erie discr\`ete); alors dans le lemme ci-dessus, on a en plus n\'ecessairement pour tout $(\rho,d,f)$, $d-f\geq 0$ (resp. $d-f>0$).
\end{rmq}
Raisonnons par l'absurde, on suppose qu'avec les notations du lemme pr\'ec\'edent, il existe un des triplets $(\rho,d,f)$ tel que $d-f<0$ et on fixe un tel triplet le plus \`a gauche possible; il ne peut \^etre en premi\`ere position car cela contredirait le crit\`ere de Casselman.  Pour tout triplet $(\rho',d',f')$ plus \`a gauche $f'\leq d'$ (par minimalit\'e) et donc n\'ecessairement $d'-f'\geq 0 >d-f$ et $f-f'>d-d'$; par les in\'egalit\'es sur les extr\'emit\'es $f$, on a aussi $f'\geq f$ d'o\`u aussi $d'\geq d$ et  le segment $[d,-f]$ est inclus dans le segment $[d',-f']$. Ainsi l'induite pour le bon groupe $GL(n,E)$:
$
\langle d',-f'\rangle_{\rho'}\times \langle d,-f\rangle_{\rho}$
est irr\'eductible et on peut donc commuter les deux facteurs. En agissant ainsi de proche en proche, on ram\`ene $\langle d,-f\rangle_{\rho}$ en premi\`ere position et on a une contradiction avec le crit\`ere de Casselman.

\subsection{Sur les points de r\'eductibilit\'e des s\'eries discr\`etes\label{reddisc}}
Soit $\pi$ une repr\'esentation temp\'er\'ee irr\'eductible; on note $Red(\pi)$ l'ensemble des suites finies de triplets $\{(\rho_{i},x_{i}, \pi_{i}), i\in [1,\ell]\}$ telle que $\ell$ est un entier, pour tout $i\in [1,\ell], x_{i}\in {\mathbb R}_{>1/2}$ et $\pi_{i}$ une repr\'esentation temp\'er\'ee irr\'eductible qui v\'erifie inductivement que $\pi_{i}\hookrightarrow \rho_{i}\vert\,\vert^{x_{i}}\times \pi_{i-1}$ avec $\pi_{0}=\pi$. On remarque que contrairement au cas des s\'eries discr\`etes strictement positives ( paragraphe \ref{seriesdiscretespositives}), on demande vraiment que $\pi_{i}$ soit un sous-module et pas seulement un sous-quotient. On demande encore que $x_{i}\geq x_{i+1}$ pour tout $i\in [1,\ell[$.

On fixe aussi une inclusion comme dans \ref{lemmedestructure}; donc on a $Jord(\pi_{+})$ et les triplets $(\rho,d,f)$, ensemble de triplets que l'on note $\mathcal T$. 

\begin{lem} Soit $\{(\rho_{i},x_{i},\pi_{i});i \in [1,\ell]\}$ un \'el\'ement de $Red(\pi)$. Alors pour tout $i\in [1,\ell]$, il existe soit $(\rho,a)\in Jord(\pi_{+})$ tel que $\rho_{i}\simeq \rho$ et $a=2x_{i}-1$ soit il existe $(\rho,d,f)\in {\mathcal T}$ tel que soit $\rho_{i}\simeq \rho$ et $x_{i}=d+1$ soit $\theta(\rho_{i})\simeq \rho$ et $x_{i}=f+1$.
\end{lem}
On peut ajouter que, pour tout $i\in [1,\ell]$, $\pi_{i}$ v\'erifie une inclusion comme dans \ref{lemmedestructure}  en rempla\c{c}ant dans celle-ci certains $(\rho,d,f)$ par $(\rho,d+1,f)$ ou $(\rho,d,f+1)$ et/ou en rempla\c{c}ant $\pi_{+}$ par une s\'erie discr\`ete totalement positive $\pi_{i,+}$ ayant comme blocs de Jordan ceux   de $\pi_{+}$ mais certains d'entre eux passant de $(\rho,a)$ \`a $(\rho,a+2)$. Comme cela et en tenant compte des in\'egalit\'e $x_{j}\geq x_{j+1}$, il suffit de faire la preuve avec $\ell=1$. On a par hypoth\`ese:
$$
\pi_{1}\hookrightarrow \rho_{1}\vert\,\vert^{x_{1}}\times \pi \hookrightarrow$$
$$
\rho_{1}\vert\,\vert^{x_{1}} \times_{(\rho,d,f)\in {\mathcal T}}\langle d,-f\rangle_{\rho} \times\pi_{+}.
$$
On d\'eplace $\rho_{1}\vert\,\vert^{x_{1}}$ vers la droite; c'est toujours possible par isomorphisme au-dessus d'une repr\'esentation $\langle d,-f\rangle_{\rho}$ si l'une au moins des deux conditions suivantes n'est pas satisfaite: $\rho_{1}\simeq \rho$ et $x_{1}=d+1$. Si les deux conditions sont satisfaites (c'est une des \'eventualit\'es du lemme) on a deux cas possibles: soit $\pi_{1}$ est inclus dans l'induite obtenue en rempla\c{c}ant $\langle d,-f\rangle_{\rho}$ par $\langle d+1,-f\rangle_{\rho}$, soit on peut quand m\^eme faire commuter $\langle d,-f\rangle_{\rho}$ et $\rho_{1}\vert\,\vert^{x_{1}}$ en gardant l'inclusion. Le premier cas, termine la preuve et dans le deuxi\`eme on continue: si on arrive jusqu'\`a $$\pi_{1}\hookrightarrow \times_{(\rho,d,f)\in {\mathcal T}}\langle d,-f\rangle_{\rho} \times \rho\vert\,\vert^{x_{1}}\times \pi_{+},
$$
on peut remplacer $\rho\vert\,\vert^{x_{1}}\times \pi_{+}$ par un sous-quotient irr\'eductible. Tous les sous-quotients irr\'eductibles sont des s\'eries discr\`etes positives sauf celui qui est inclus dans 
$\theta(\rho_{1})\vert\,\vert^{-x_{1}}\times \pi_{+}$. Sauf le dernier cas, on est encore dans l'une des \'eventualit\'es du lemme gr\^ace \`a \ref{seriesdiscretespositives}. Reste le dernier cas, o\`u on pousse $\theta(\rho_{1})\vert\,\vert^{-x_{1}}$ vers la gauche; on obtient soit l'un des cas restant du lemme, soit $$\pi_{1}\hookrightarrow \theta(\rho)_{1}\vert\,\vert^{-x_{1}}\times _{(\rho,d,f)\in {\mathcal T}}\langle d,-f\rangle_{\rho} \times\pi_{+}.
$$
Ceci est impossible car $\pi_{1}$ est suppos\'ee une repr\'esentation temp\'er\'ee.
\begin{cor} Soit $\pi$ une repr\'esentation temp\'er\'ee irr\'eductible et soit $\{(\rho_{i},x_{i},\pi_{i});$ $i \in [1,\ell]\}$ un \'el\'ement de $Red(\pi)$. Alors $\sum_{i\in [1,\ell]}d_{\rho_{i}}(2x_{i}-1)\leq n$; supposons que l'on ait  l'\'egalit\'e et  que pour tout $i\in [1,\ell]$,  $\rho_{i}\simeq \theta(\rho_{i})$ et $x_{i}$ est un demi-entier, alors le support cuspidal \'etendu de $\pi$ est le support cuspidal de la repr\'esentation temp\'er\'ee $\times_{i\in [1,\ell]}St(\rho_{i},2x_{i}-1)$ de $GL(n,E)$.
\end{cor}
En tenant compte du lemme pr\'ec\'edent (cf. aussi la preuve de ce lemme), on a
$$
\sum_{i\in [1,\ell]}d_{\rho_{i}}(2x_{i}-1)\leq \sum_{(\rho,d,f)\in {\mathcal T}} d_{\rho}(2d+1+2f+1)+\sum_{(\rho,a)\in Jord(\pi_{+})}d_{\rho}a.$$
Evidemment $2d+1+2f+1=2(d+f+1)$ vaut 2 fois la longueur du segment $[d,-f]$ et le terme de droite vaut donc exactement $n$. L'\'egalit\'e force le fait que pour tout $(\rho,d,f)\in {\mathcal T}$, il existe $i\in [1,\ell]$ tel que $\rho_{i}\simeq \rho$ et $x_{i}=d+1$ et il existe $j\in [1,\ell]$ tel que $\rho_{j}\simeq \theta(\rho_{j})=\rho$ et $x_{j}=f+1$ et que pour tout $(\rho,a)\in Jord(\pi_{+})$ il existe $i'\in [1,\ell]$ tel que $\rho_{i'}\simeq \rho$ et $2x_{i'}-1=a$. Ce qui nous int\'eresse, avec les hypoth\`eses de la fin du lemme est que cela force tout triplet $(\rho,d,f)\in {\mathcal T}$ de v\'erifier $\rho\simeq \theta(\rho)$ et $d,f\in 1/2 {\mathbb Z}_{\geq 0}$ puisque le suppose que les $\rho_{i}\simeq \theta(\rho_{i})$ et $x_{i}\in 1/2{\mathbb N}$ et $x_{i}>1/2$.  Par d\'efinition du support cuspidal \'etendu, celui de $\pi$ est l'union pour $(\rho,d,f)\in {\mathcal T}$ des ensembles $\rho\vert\,\vert^x; x\in [d,-f]$ et $\theta(\rho)\vert\,\vert^{-x}; x\in [d,-f]$ auquel on ajoute le support cuspidal \'etendu de $\pi_{+}$.  On peut donc r\'ecrire les ensembles ci-dessus sous la forme
$
\rho\vert\,\vert^x$ avec $x\in [d,-d]\cup [f,-f]$. Ceci est le support cuspidal de la repr\'esentation temp\'er\'ee $St(\rho,2d+1)\times St(\rho,2f+1)$. Comme on sait d\'ej\`a que le support cuspidal \'etendu de $\pi_{+}$ est celui d'une repr\'esentation temp\'er\'ee, cela conclut la preuve du lemme.

\subsection{Paquet stable de s\'eries discr\`etes\label{paquetstableseriesdiscretes}}
\begin{thm}
Soit $\pi$ une s\'erie discr\`ete de $G$; soit aussi $\pi^{GL}$ une repr\'esentation $\theta$-elliptique  de $\tilde{GL}$  dont la $\theta$-trace a une projection sur $I_{cusp}^{\underline{G},st}$ qui n'est pas ortho\-gonale \`a  $\pi$. Alors le support cuspidal de $\pi^{GL}$ co\"{\i}ncide avec le support cuspidal \'etendu de $\pi$.
\end{thm}
On va d\'emontrer en m\^eme temps le compl\'ement important suivant:
\begin{prop}
Soit $\pi$ une repr\'esentation elliptique de $G$ dont la projection de la trace sur $I_{cusp}^{G,st}$ n'est pas nulle alors $\pi$ est une s\'erie discr\`ete.
\end{prop} On proc\`ede ainsi: on fixe une repr\'esentation $\theta$-elliptique de $\tilde{GL}(n,E)$ que l'on note $\pi^{GL}$; on suppose que la $\theta$-trace de $\pi^{GL}$ a une projection non nulle sur $I_{cusp}^{G,st}$ et on montre que toute repr\'esentation elliptique intervenant dans cette $\theta$ trace est une s\'erie discr\`ete et que le support cuspidal \'etendu de cette s\'erie discr\`ete est le support cuspidal de $\pi^{GL}$.

On \'ecrit $\pi^{GL} $ sous la forme $\times_{(\rho_{i},a_{i})\in i\in [1\ell]}St(\rho_{i},a_{i})$ et on suppose que pour $i\in [1,\ell[$, $a_{i}\leq a_{i+1}$. On fixe $\pi$ une repr\'esentation temp\'er\'ee irr\'eductible qui compose une des repr\'esentations elliptiques intervenant dans la projection sur $I_{cusp, st}^G$ de cette $\theta$-trace. On va montrer que pour tout $i\in [1,\ell]$, il existe une repr\'esentation temp\'er\'ee irr\'eductible $\pi_{i}$ tel que $\{(\rho_{i},(a_{i}+1)/2,\pi_{i});i\in [1,\ell]$ soit dans $Red(\pi)$. V\'erifions que cela suffira: on sait par construction que 
$$
\sum_{i\in [1,\ell]} d_{\rho_{i}}a_{i}=n.
$$
Ainsi le support cuspidal \'etendu de $\pi$ est le support cuspidal de $\pi^{GL}$; comme $\pi^{GL}$ est elliptique tous les $(\rho_{i},a_{i})$ sont distincts avec $\rho_{i}\simeq \theta(\rho_{i})$ pour tout $i$; et  $\pi$ est n\'ecessairement une s\'erie discr\`ete: en effet  une repr\'esentation temp\'er\'ee qui n'est pas une s\'erie discr\`ete  v\'erifie une inclusion comme en \ref{lemmedestructure} (1) mais avec un des triplets $(\rho,d,f)$ v\'erifiants $d=f$; cela force l'existence de $i\neq i'\in [1,\ell]$ avec $\rho_{i}\simeq \rho$ et $a_{i}=2d-1$ et $\theta(\rho_{i'})=\rho$ avec $a_{i'}=f+1$, d'o\`u $\rho_{i}=\rho_{i'}$ et $a_{i}=a_{i'}$ ce qui a \'et\'e exclu.

Il reste donc \`a construire les $\pi_{i}$ pour $i\in [1,\ell]$; on proc\`ede de fa\c{c}on descendante. On pose $n_{+}:= n+\sum_{i\in [1,\ell]}2d_{\rho_{i}}$ et on note $\pi^{GL}_{+}$ la repr\'esentation de $GL(n_{+},E)$ \'egale \`a $\times_{i\in [1,\ell]}St(\rho_{i},a_{i}+2)$. On consid\`ere la projection de la $\theta$-trace de $\pi^{GL}_{+}$ sur $I_{cusp}^{G_{n_{+}},st}$; on sait d'apr\`es \ref{transfertetjac} qu'en appliquant $$\circ_{i\in [\ell,1]}Jac_{\rho_{i}\vert\,\vert^{(a_{i}+1)/2}}$$ \`a cette projection, on trouve la projection de la $\theta$-trace de $\pi^{GL}$ sur $I_{cusp}^{G_{n},st}$. En tenant maintenant compte de \ref{pasdeuxfois}, on voit qu'il existe une  repr\'esentation not\'ee $\pi_{\ell}$ temp\'er\'ee irr\'eductible qui est composante d'une repr\'esentation elliptique intervenant dans la projection de $\pi^{GL}_{+}$ sur $I_{cusp}^{G}$ et telle que $$\circ_{i\in [\ell,1]}Jac_{\rho_{i}\vert\,\vert^{(a_{i}+1)/2}}\pi_{\ell}=\pi.$$ On pose, pour tout $i\in [1,\ell[$,
$
\pi_{i}=\circ_{j\in [\ell,i+1]}Jac_{\rho_{j}\vert\,\vert^{(a_{j}+1)/2}}\pi_{\ell}.
$
La repr\'esentation $\pi_{i}$ est n\'ecessairement temp\'er\'ee car elle intervient dans la projection de $$\times_{j\in [1,i]} St(\rho_{i},a_{i}+2)\times_{j\in [i+1,\ell]}St(\rho_{i},a_{i})$$  sur $I_{cusp}^{G_{n+2\sum_{i\in [1,i]}d_{\rho_{i}}},st}$ et on a inductivement
$$
Jac_{\rho_{i+1}\vert\,\vert^{(a_{i}+1)/2}}\pi_{i+1}=\pi_{i}.
$$Ainsi $\pi_{i+1}$ est un sous-module irr\'eductible de l'induite $\rho_{i+1}\vert\,\vert^{(a_{i+1}+1)/2}\times \pi_{i}$.
Avec \ref{pasdeuxfois}, on sait que $$Jac_{\rho_{i+1}\vert\,\vert^{(a_{i+1}+1)/2},\rho_{i+1}\vert\,\vert^{(a_{i+1}+1)/2}}\pi_{i+1}=0$$ et donc que $Jac_{\rho_{i+1}\vert\,\vert^{(a_{i+1}+1)/2}}\pi_{i}=0$. Ainsi par r\'eciprocit\'e de Frobenius  l'induite $\rho_{i+1}\vert\,\vert^{(a_{i+1}+1)/2}\times \pi_{i}$ a un unique sous-module irr\'eductible qui est n\'ecessairement $\pi_{i+1}$. On a donc toutes les propri\'et\'es souhait\'ees.

 \begin{cor} Soit $\pi$ une s\'erie discr\`ete de $G$. 
Il existe une unique  combinaison lin\'eaire de toutes les s\'eries discr\`etes de $G$ ayant m\^eme support cuspidal \'etendu que $\pi$ qui soit stable  et cette combinaison lin\'eaire se transf\`ere en l'unique repr\'esentation $\theta$-elliptique de $\tilde{GL}(n,E)$ ayant pour support cuspidal le support cuspidal \'etendu de $\pi$. De plus $I_{cusp}^{G,st}$ est engendr\'e comme espace vectoriel par ces combinaisons lin\'eaires stables.
\end{cor}
On fixe $\pi$; on a donc le support cuspidal \'etendu de $\pi$ d'o\`u l'existence d'une unique repr\'esentation $\theta$-elliptique de $\tilde{GL}$ ayant pour support cuspidal le support cuspidal \'etendu de $\pi$; on la note $\pi ^{GL}$. On sait que $\pi^{GL}$ est un transfert d'une repr\'esentation virtuelle stable $\pi_{st}$ de $G_{n}$; $\pi_{st}$ est orthogonal \`a toute s\'erie discr\`ete n'ayant pas m\^eme support cuspidal \'etendu que $\pi$ d'apr\`es le th\'eor\`eme pr\'ec\'edent; cette combinaison lin\'eaire  ne peut \^etre orthogonale \`a une s\'erie discr\`ete, $\pi'$, ayant m\^eme support cuspidal \'etendu que $\pi$ car sinon $\pi'$ aurait une projection nulle sur $I_{cusp}^{\underline{G},st}$ ce qui a \'et\'e exclu en \ref{projectionnonnulle}.

\subsection{Support cuspidal \'etendu et endoscopie \label{supportcuspidaletendoscopie}}
\begin{thm} 
Soit $\underline{H}$ une donn\'ee endoscopique elliptique de $G$ et $\Pi_{\underline{H},st}$ une combinaison lin\'eaire  stable de s\'eries discr\`etes de $\underline{H}$ ayant toutes m\^eme support cuspidal \'etendu; on suppose que cette combinaison lin\'eaire se transforme sous les automorphismes du groupe endoscopique suivant le caract\`ere d\'etermin\'e par $G$ (le caract\`ere trivial, ici). Le transfert endoscopique de cette distribution stable  est une combinaison lin\'eaire de repr\'esentations elliptiques de $G$ dont toutes les composantes irr\'eductibles ont pour support cuspidal \'etendu l'image (cf. la preuve) du support cuspidal \'etendu des repr\'esentations de $\underline{H}$ dont on est parti.
 \end{thm}

Le lemme \ref{pasdeuxfois} n'est pas exact pour une repr\'esentation elliptique; dans cet \'enonc\'e on prenait comme hypoth\`ese que la repr\'esentation elliptique appara\^{\i}t dans une combinaison lin\'eaire stable de repr\'esentations elliptiques; or nous avons vu en \ref{paquetstableseriesdiscretes} que cette hypoth\`ese n'est pas satisfaite en g\'en\'eral. Il faut \`a la place utiliser le lemme suivant, pour \'eviter des annulations dans le calcul de modules de Jacquet:
\begin{lem}Soit $\pi$ une repr\'esentation temp\'er\'ee irr\'eductible composante d'une re\-pr\'e\-sentation elliptique de $G$. Et soit $\rho$ une repr\'esentation cuspidale unitaire de $GL(d_{\rho},E)$ et $x\in {\mathbb R}$. Alors $Jac_{\rho\vert\,\vert^x}\circ Jac_{\rho\vert\,\vert^x}\pi$ est soit nul soit est une repr\'esentation  irr\'eductible, $\pi'$, et on a alors que $\pi$ est l'unique sous-repr\'esenta\-tion irr\'eductible de l'induite $\rho\vert\,\vert^x\times \rho\vert\,\vert^x\times \pi'$. 
\end{lem}
Il existe n\'ecessairement une donn\'ee endoscopique elliptique $\underline{H}$ de $G$ (qui peut \^etre $G$ lui-m\^eme) et une combinaison lin\'eaire stables de s\'eries discr\`etes de $\underline{H}$, $\Pi_{\underline{H},st}$, vu comme un \'el\'ement de $I_{cusp}^{\underline{H},st}$, se transformant correctement sous les automorphismes de la donn\'ee endoscopique, dont le transfert \`a $G$ a $\pi$ comme l'une de ses composantes.

En proc\'edant comme dans la preuve de \ref{pasdeuxfois}, on note $t$ le plus grand entier tel qu'il existe $\tau$ une composante du transfert tel que $Jac_{\rho\vert\,\vert^x}\circ \cdots \circ Jac_{\rho\vert\,\vert^x}\tau$ est non nul o\`u il y a $t$ facteurs $Jac_{\rho\vert\,\vert^x}$. On v\'erifie ici que $t$ est inf\'erieur ou \'egal au plus grand entier $t'$ tel que $Jac_{\rho\vert\,\vert^x}\circ \cdots \circ Jac_{\rho\vert\,\vert^x}\Pi_{\underline{H},st}$ soit non nul, o\`u il y a $t'$ facteur $Jac_{\rho\vert\,\vert^x}$. On calcule $t'$ en utilisant \ref{pasdeuxfois}; $t'$ vaut au plus 2, puisque $\underline{H}$ est un produit d'au plus deux groupes. D'o\`u $t\leq 2$. Le cas $t=2$ est tout \`a fait possible. La suite du lemme est comme en \ref{pasdeuxfois}. 

\

On d\'efinit ici $Red(\pi)$ comme \'etant l'ensemble des suites, de la forme $$(\rho_{i},x_{i},\pi_{i}); i\in [1,\ell]$$ ($\ell$ d\'epend de la suite) tel que pour tout $i\in [1,\ell]$ $x_{i}\geq x_{i+1}>1/2$, $\pi_{i}$ est une repr\'esentation temp\'er\'ee (on pose $\pi_{0}=\pi$) telle que $\pi_{i+1}=\pi_{i}$ si $(\rho_{i},x_{i})=(\rho_{i+1},x_{i+1})$ et $\pi_{i}\hookrightarrow \rho_{i}\vert\,\vert^{x_{i}}\times \pi_{i-1}$ si $(\rho_{i},x_{i})\neq (\rho_{i+1},x_{i+1})$ et dans le cas restant:
$$
\pi_{i+1}=\pi_{i}\hookrightarrow \rho_{i}\vert\,\vert^{x_{i}}\times  \rho_{i}\vert\,\vert^{x_{i}}\times \pi_{i-1}.$$
Ensuite on d\'emontre comme dans \ref{reddisc} que pour tout \'el\'ement comme ci-dessus de $Red(\pi)$, on a $\sum_{i\in [1,\ell]}d_{\rho_{i}}(2x_{i}-1)\leq n$. La fin du corollaire de loc. cite est aussi exacte, c'est-\`a-dire que  l'\'egalit\'e coupl\'ee avec le fait que les $x_{i}$ sont des demi-entiers et $\rho_{i}\simeq \theta(\rho_{i})$ entra\^{\i}nent que le support cuspidal de $\pi$ est le support cuspidal de la repr\'esentation temp\'er\'ee
$\times_{i\in [1,\ell]}St(\rho_{i},2x_{i}-1)$.

On d\'emontre maintenant le th\'eor\`eme; on va le sp\'ecialiser au cas du groupe unitaire car c'est dans ce cas que l'image des groupes endoscopiques fait intervenir quelques torsions. On fixe donc une d\'ecomposition $n=n_{1}+n_{2}$ et $\underline{H}$ la donn\'ee endoscopique associ\'e (cf\cite{inventaire} 1.8); pour nous, on voit $\underline{H}$ comme un produit de groupes unitaires, $H:=U(n_{1},E/F)\times U(n_{2},E/F)$ et la donn\'ee de deux  caract\`eres de $E^*$, $\omega_{i}$ pour $i=1,2$ tel que pour $i=1,2$,  la restriction de $\omega_{i}$ \`a $F^*$ vaut le caract\`ere quadratique associ\'e \`a l'extension $E$ de $F$ si $n-n_{i}$ est impaire et est le caract\`ere trivial sinon. On suppose que $\omega_{1}=\omega_{2}$ si $n_{1}=n_{2}$. Il y a un automorphisme non trivial de la donn\'ee endoscopique que si $n_{1}=n_{2}$ et l'automorphisme non trivial \'echange de fa\c{c}on \'evidente les deux facteurs.

Un paquet stable de s\'eries discr\`etes, pour $U(n_{1},E/F)\times U(n_{2},E/F)$ est associ\'e \`a certaines repr\'esentations temp\'er\'ees de $GL(n_{1},E)\times GL(n_{2},F)$ que l'on a donc le droit d'\'ecrire sous la forme
$$
\times_{(\rho,a)\in {\mathcal E}_{1}}St(\rho\otimes \omega_{1},a)\times \times_{(\rho,a)\in {\mathcal E}_{2}}St(\rho\otimes \omega_{2},a). \eqno(1)
$$
En ayant ainsi gliss\'e les torsions dans la d\'efinition de ${\mathcal E}_{i}$ pour $i=1,2$ l'image du support cuspidal pour $GL(n,E)$ est tout simplement le support cuspidal de la repr\'esentation 
$\times_{(\rho,a)\in {\mathcal E}_{1}\cup {\mathcal E}_{2}}St(\rho,a)$.

Le paquet (1) n'est stable sous l'automorphisme non trivial quand il existe que si ${\mathcal E}_{1}={\mathcal E}_{2}$ et si cette condition est v\'erifi\'ee, la distribution stable associ\'ee au paquet de s\'eries discr\`etes de $U(n/2,E/F)\times U(n/2,E/F)$ ayant comme support cuspidal \'etendu le support cuspidal de (1) est n\'ecessairement invariant par l'automorphisme par unicit\'e. Dans les autres cas, il faut \'eventuellement sommer deux paquets de repr\'esentations.

Le transfert entrelace $Jac_{\rho\vert\,\vert^x}$ avec $Jac_{\rho\otimes \omega_{1}\vert\,\vert^x}\otimes Id \oplus Id\otimes Jac_{\rho\otimes \omega_{2}\vert\,\vert^x}$. C'est l\`a que l'on voit que comme ${\mathcal E}_{1}\cup {\mathcal E}_{2}$ peut avoir de la multiplicit\'e 2, on peut tr\`es bien avoir $t'=2$ avec les notations du d\'ebut de la preuve.

On fixe ${\mathcal E}_{i}$ pour $i=1,2$ comme dans (1) et on d\'emontre le th\'eor\`eme pour la distribution stable associ\'ee rendue, si n\'ecessaire, invariante sous le groupe d'automorphisme de la donn\'ee endoscopique. On pose ${\mathcal E}={\mathcal E}_{1}\cup {\mathcal E}_{2}$ et on ordonne ${\mathcal E}$ en notant $\ell$ le nombre d'\'el\'ements de ${\mathcal E}$ en prenant en compte les multiplicit\'es \'eventuelles et en \'ecrivant:
$$
{\mathcal E}=\{(\rho_{i},a_{i}); i\in [1,\ell]\},
$$
l'ordre v\'efiiant: s'il existe $i\neq j\in [1,\ell]$ tel que $\rho_{i}=\rho_{j}$ et $a_{i}=a_{j}$ alors $\vert i-j\vert=1$ et  $a_{1}\geq \cdots \geq a_{\ell}$.

On consid\`ere $n_{+}=n+2\sum_{i\in [1,\ell]}d_{\rho_{i}}$ et pour $i=1,2$, $$n_{i,+}=n_{i}+2\sum_{(\rho,a)\in {\mathcal E}_{i}}d_{\rho}.$$ Pour la donn\'ee endoscopique $U(n_{1,+},E/F)\times U(n_{2,+},E/F)$ de $U(n_{+},E/F)$ on consid\`ere la distribution stable combinaisons lin\'eaires de s\'eries discr\`etes  associ\'ee (par le corollaire \ref{paquetstableseriesdiscretes}) \`a ${\mathcal E}_{1,+}\times {\mathcal E}_{2,+}$ o\`u pour $i=1,2$, on obtient ${\mathcal E}_{i,+}$ en rempla\c{c}ant $(\rho,a)$ par $(\rho,a+2)$ dans ${\mathcal E}_{i}$. On rend \'eventuellement invariant ce paquet par l'automorphisme de la donn\'ee endoscopique et on le transf\`ere; la situation de d\'epart s'obtient \`a partir de cette situation en appliquant $\circ_{i\in [1,\ell]}Jac_{\rho_{i}\vert\,\vert^{(a_{i}+1)/2}}$ \`a un scalaire pr\`es (qui d\'epend des multiplicit\'es) mais qui n'importe pas. Ainsi, en tenant compte du lemme ci-dessus, dans ce transfert, il existe une composante $\pi_{+}$ tel qu'\`a un scalaire pr\`es:
$$
\pi=\circ_{i\in [1,\ell]}Jac_{\rho_{i}\vert\,\vert^{(a_{i}+1)/2}}\pi_{+}.
$$
Ensuite on termine la d\'emonstration comme dans le cas des s\'eries discr\`etes.
\subsection{D\'efinition des paquets de Langlands des s\'eries dis\-cr\`etes\label{deflang}}
\begin{thm}
Le support cuspidal \'etendu est l'invariant d\'eterminant les paquets de Langlands. Ou encore, la combinaison lin\'eaire stable de s\'erie discr\`ete d\'ecrite dans le corollaire \ref{paquetstableseriesdiscretes} est la projection de la trace de $\pi$ sur $I_{cusp}^{\underline{G},st}$.
\end{thm}
\begin{rmq}
Soit $\pi$ une repr\'esentation elliptique de $G$ et soit ${\underline{H}}$ une donn\'ee endoscopique elliptique $G$; alors  la projection de la trace de $\pi$ sur $I_{cusp}^{\underline{H},st}$ est soit 0 soit exactement (\`a un scalaire pr\`es) le paquet invariant sous le groupe d'automorphisme de $\underline{H}$ port\'e par les s\'eries discr\`etes de $\underline{H}$ dont l'image du support cuspidal \'etendu est exactement le support cuspidal \'etendu de $\pi$.
\end{rmq}
La remarque est plus g\'en\'erale que le th\'eor\`eme c'est donc elle que nous allons montrer mais en fait elle n'apporte pas grand chose car  dans l'\'enonc\'e on peut tr\`es certainement enlever le ''est soit 0''  et ajouter \`a la fin quand un tel paquet existe. On ne d\'emontre pas ce r\'esultat plus g\'en\'eral ici car la d\'emonstration potentielle doit imiter celle d'Arthur et en m\^eme temps calculer les coefficients de ces projections; cela n\'ecessite de repasser \`a une situation globale et donc d\'epasse le cadre de ce travail; le cas des groupes orthogonaux et symplectiques est fait par \cite{book} et celui des groupes unitaires est dans \cite{white} et \cite{mok}. 

On fixe une repr\'esentation elliptique de $G$. Pour toute donn\'ee endoscopique elliptique $\underline{H}$ (prise \`a \'equivalence pr\`es) de $G$ on note $\pi^{H,st}$ la projection de la trace de $\pi$ sur $I_{cusp}^{\underline{H},st}$ et on d\'ecompose $\pi^{H,st}$ en $\pi^{H,st}_{0}+\pi^{H,st}_{1}$, o\`u $\pi^{H,st}_{1}$ est la projection de $\pi^{H,st}$ sur l'espace orthogonal \`a toute trace de s\'erie discr\`ete de $\underline{H}$ dont le support cuspidal \'etendu n'est pas celui de $\pi$; ainsi $\pi^{H,st}_{i}$ est  invariant par le groupe d'automorphisme de $\underline{H}$ pour $i=0,1$. On sait avec \ref{supportcuspidaletendoscopie} que le transfert de $\oplus_{\underline{H}}\pi^{H,st}_{1}$ est un \'el\'ement de $I_{cusp}^{G}$ orthogonal \`a la trace de $\pi$. Donc la trace de $\pi$ est le transfert de $\oplus_{\underline{H}}\pi^{H,st}_{0}$ et le r\'esultat.
\subsection{Conclusion\label{conclusion}}
Etant donn\'e une s\'erie discr\`ete $\pi$ de $G$, on a associ\'e \`a $\pi$ deux distributions stables: l'une est obtenue naturellement en projetant un pseudocoefficient cuspidal de $\pi$ dans $I_{cusp}^{G,st}$ (\ref{pseudocoefficient}) et l'autre en consid\'erant l'unique \'el\'ement de $I_{cusp}^{G,st}$ non orthogonal \`a $\pi$ et dont le transfert \`a $I_{cusp}^{\tilde{GL}}$ est une repr\'esentation irr\'eductible. Le th\'eor\`eme ci-dessus montre que ces deux d\'efinitions co\"{\i}ncident. On d\'efinit ainsi le paquet de Langlands de $\pi$ comme l'ensemble des s\'eries discr\`etes $\pi'$ de $G$ qui constituent la distribution stable pr\'ec\'edente; ce sont pr\'ecis\'ement les s\'eries discr\`etes ayant le m\^eme support cuspidal \'etendu que $\pi$. Ainsi tous les paquets de Langlands sont disjoints.

On \'etend ces d\'efinitions aux cas des repr\'esentations temp\'er\'ees en utilisant le fait que l'induction commute au transfert.

\section{Morphisme dans le L-groupe}
\noindent
Maintenant que l'on a d\'emontr\'e que les paquets de Langlands des s\'eries discr\`etes de $G$ sont param\'etr\'es par certaines repr\'esentations temp\'er\'ees de $\tilde{GL}(n,E)$, on peut associer \`a un tel paquet le morphisme de Langlands de $W_{F}\times SL(2,{\mathbb C})$ dans le $L$-groupe de ${GL}(n,E)$ vu comme groupe d\'efini sur $F$ (si $E\neq F$). Le but de cette partie est de d\'emontrer que l'image, \`a conjugaison pr\`es, se factorise par le $L$-groupe de $G$. Cet aspect est fait dans \cite{book} en quelques pages et il est tout \`a fait exact qu'il n'y a rien de profond dans la d\'emonstration; toutefois la d\'emonstration utilise elle un r\'esultat profond qui est l'\'egalit\'e des fonctions $L$ dites d'Artin (qui ont \'et\'e d\'efinies par Shahidi) avec les fonctions $L$ que l'on obtient via la correspondance de Langlands. Ce r\'esultat est d\^u \`a H\'enniart en \cite{henniartasai} et utilise les r\'esultats de Shahidi \cite{shahidiannals}.

\subsection{Le cas des morphismes irr\'eductibles\label{morphismeirr}}
\begin{prop} Soit $\pi^{GL}$ une s\'erie discr\`ete irr\'eductible de $\tilde{GL}$ qui est $\theta$ invariante. Alors il existe un unique sous-groupe endoscopique elliptique simple  tel que le param\`etre de Langlands de $\pi^{GL}$ se factorise, \`a conjugaison pr\`es, par son $L$-groupe
\end{prop}
Le cas de $GL(n,F)$ et de son automorphisme ext\'erieur est bien connu: un tel morphisme est soit orthogonal, soit symplectique, cela se voit sur l'existence ou non de p\^ole \`a la fonction $L(\pi^{GL},r,s)$ en $s=0$ pour $r$ soit la repr\'esenation $Sym^2 {\mathbb C}^n$ soit la repr\'esentaiton $\wedge^2 {\mathbb C}^n$  de $GL(n^*,{\mathbb C})$; un p\^ole \`a la fonction $L$ associ\'e \`a $Sym^2$ ne suffit pas pour d\'eterminer le $L$ groupe, ce qui est bien connu depuis longtemps: l'existence du p\^ole assure que le morphisme se factorise par $O(n,{\mathbb C})$. En composant avec le d\'eterminant, on obtient un caract\`ere quadratique qui, si $n$ est pair, d\'etermine le discriminant du groupe orthogonal (et donc l'action du groupe de Galois, pour le $L$-groupe) et  si $n$ est impair est le caract\`ere $\chi$ de \cite{inventaire} 1.8.   Cela se g\'en\'eralise  (cf. \cite{gsp4} paragraphe 2, page 68) au cas de $GL(n,F)\times GL(1,F)$.

Le cas moins connu est celui o\`u $E$ est une extension quadratique de $F$ donnant lieu aux fonctions $L$ d'Asai. Et on va le traiter ici mais  le r\'esultat n'est \'evidemment pas  nouveau, bien au contraire (cf. par exemple \cite{rogawski},\cite{bellaichechenevier}).

\subsection{Le cas des groupes unitaires et de leur groupe dual\label{casunitaire}}
Pour la commodit\'e du lecteur on r\'ecrit ici le formalisme de Langlands concernant les groupes unitaires.

\

Soit $\pi$ une repr\'esentation de $GL(n,E)$; on lui associe, via la fonctorialit\'e de Langlands, un morphisme de $W_{E}\times SL(2,{\mathbb C})$ dans $GL(n,{\mathbb C})$. On note $\psi_{\pi,E}$ le morphisme ainsi associ\'e.

\begin{rmq} On suppose que $\pi$ est une s\'erie discr\`ete autocontragr\'ediente; alors il existe un signe $\lambda_{\pi}$ tel que pour tout $\sigma\in W_{F}-W_{E}$, il existe $\alpha\in GL(n,{\mathbb C})$ tel que 
$$
\psi_{\pi,E}(\sigma^2)=\lambda_{\pi} \, ^t\alpha^{-1}\, \alpha\eqno(1)$$
$$
\forall w\in W_{E}\times SL(2,{\mathbb C}), \psi_{\pi,E}(\sigma^{-1} w\sigma)=\alpha^{-1} \, ^t\psi_{\pi,E}(w)^{-1}\alpha.\eqno(2)
$$
\end{rmq}
On  v\'erifie d'abord que si $\lambda_{\pi}$ existe tel que (1) et (2) sont v\'erifi\'es pour un choix de $\sigma$ alors, pour ce $\lambda_{\pi}$, (1) et (2) sont v\'erifi\'es pour tout choix de $\sigma$, bien s\^ur $\alpha$ d\'epend lui de $\sigma$. L'unicit\'e de $\lambda_{\pi}$ est aussi claire car par hypoth\`ese $\pi$ est une s\'erie discr\`ete et $\psi_{\pi,E}$ d\'efinit donc une repr\'esentation irr\'eductible de $W_{E}$; ainsi $\alpha$ en (2) est d\'etermin\'e \`a un scalaire pr\`es.

Montrons l'existence de $\lambda_{\pi}$.
On fixe donc $\sigma$ et $\alpha$ tel que (2) soit v\'erifi\'e, ce qui est possible puisque $\pi$ est suppos\'e autocontragr\'ediente. On pose $\beta:= \psi_{\pi,E}(\sigma^2)$. En rempla\c{c}ant dans (2), $w$ par $\sigma^{-1}w \sigma$, on obtient, pour tout $w\in W_{E}\times SL(2,{\mathbb C})$:
$$
\beta^{-1} \psi_{\pi,E}(w)\beta=\alpha^{-1} \, ^t\psi_{\pi,E}(\sigma^{-1}w\sigma)^{-1}\alpha$$
et en appliquant encore (2), cela vaut:
$
\alpha^{-1} \, ^t\alpha \psi_{\pi,E}(w) \, ^t\alpha^{-1}\alpha.
$
Par irr\'eductibilit\'e de $\psi_{\pi,E}$, il existe donc un scalaire $\lambda$ tel que
$$
\beta=\lambda\, ^t\alpha^{-1}\alpha.
$$
On applique (2) \`a $w=\sigma^{2}$ et on obtient
$
\beta= \alpha^{-1} \, ^t \beta^{-1} \alpha$ et en rempla\c{c}ant $\beta$ par la valeur d\'ej\`a trouv\'ee:
$$
\lambda \, ^t\alpha^{-1}\alpha= \lambda^{-1} \, \alpha^{-1} \alpha \, ^t\alpha^{-1}\alpha;
$$
c'est-\`a-dire $\lambda^{2}=1$. Et $\lambda$ satisfait (1).
\begin{defi} Soit $\pi$ une repr\'esentation  $\theta$-elliptique de $GL(n,E)$; ainsi   $\pi$ est une induite de s\'eries discr\`etes, chaque s\'erie discr\`ete \'etant elle-m\^eme autocontragr\'ediente et n'intervenant qu'une fois. On parle donc du signe attach\'e \`a chaque s\'erie discr\`ete composant $\pi$.
\end{defi}
\begin{prop} Soit $\pi$ une repr\'esentation  $\theta$ elliptique de $GL(n,E)$. Alors le param\`etre de Langlands de $\pi$ se factorise, \`a conjugaison pr\`es, par le $L$ groupe de $U(n,E/F)$ si et seulement si le signe associ\'e \`a chaque s\'erie discr\`ete composant $\pi$ vaut $(-1)^{n-1}$.
\end{prop}
On note $\sigma_{0}$ l'\'el\'ement non trivial du groupe de Galois de $E/F$ et $\theta^*$ l'automor\-phis\-me dual de $\theta$.

On commence par rappeler que le $L$-groupe de $GL(n,E)$ vu comme $F$-groupe est isomorphe \`a 
$$
(GL(n,{\mathbb C})\times GL(n,{\mathbb C}))\rtimes {\rm Gal}(E/F),\eqno(1)
$$
o\`u $\sigma_{0}$ agit en permutant les deux facteurs. L'action de $\theta^*$ est alors donn\'ee par:
$$
\theta^*(g,g')=(J\, ^tg^{'-1}J^{-1}, J\, ^tg^{-1}J^{-1}),
$$
o\`u $J$ est la matrice antidiagonale ayant alternativement des 1 et des -1; en particulier $J^{-1}=^tJ= (-1)^{n-1}J$.

Le groupe dual de $U(n,E/F)$ est exactement l'ensemble des \'el\'ements invariants par $\theta^*$; il est isomorphe \`a $GL(n,{\mathbb C})\rtimes {\rm Gal}(E/F)$ o\`u $\sigma_{0}$ agit sur $GL(n,{\mathbb C})$ par $\sigma_{0}(g)=J\, ^tg^{-1}J^{-1}$, par l'isomorphisme, valant l'identit\'e sur $Gal(E/F)$ et v\'erifiant
$$
g\in GL(n,{\mathbb C}) \mapsto (g, J ^tg^{-1}J^{-1}).$$

On associe \`a $\pi$ son param\`etre de Langlands $\psi_{\pi,E}$. On suppose d'abord que toute s\'erie discr\`ete constituant $\pi$ a pour signe $(-1)^{n-1}$ et on va montrer que le param\`etre de Langlands de $\pi$ se factorise, \`a conjugaison pr\`es, par le $L$-groupe de $U(n,E/F)$.

On note $\psi_{\pi,E}$ le morphisme de Langlands associ\'e \`a $\pi$ vu comme morphisme de $W_{E}\times SL(2,{\mathbb C})$ dans $GL(n,{\mathbb C})$. Pour \'ecrire le morphisme de Langlands de $\pi$ quand on voit $GL(n,E)$ comme un groupe sur $F$, on fixe $\sigma$ un \'el\'ement de $W_{E}-W_{F}$ et on note $\psi_{\pi,F}$ le morphisme de $W_{F}\times SL(2,{\mathbb C})$ dans (1):
$$
\psi_{\pi,F}(w)=(\psi_{\pi,E}(w),\psi_{\pi,E}(\sigma^{-1}w \sigma)),$$
$$
\psi_{\pi,F}(\sigma)= \sigma_{0} (1,\psi_{\pi,E}(\sigma^2)).
$$
Pour \'eviter des erreurs on fait les v\'erifications usuelles: 
$$\psi_{\pi,F}(\sigma)^{-1}\psi_{\pi,F}(w)\psi_{\pi,F}(\sigma)=
(1,\psi_{\pi,E}(\sigma^{-2}))(\psi(\sigma^{-1}w\sigma),\psi_{\pi,E}(w))(1,\psi_{\pi,E}(\sigma^2))$$
$$
=(\psi_{\pi,E}(\sigma w \sigma^{-1}),\psi_{\pi,E}(\sigma^{-2}w\sigma^{2}))=\psi_{\pi,F}(\sigma^{-1} w \sigma);
$$
$$
\psi_{\pi,F}(\sigma)^{2}=(\psi_{\pi,E}(\sigma^2),\psi_{\pi,E}(\sigma^2))=\psi_{\pi,F}(\sigma^2).
$$
De plus le choix de $\sigma$ n'influe qu'\`a conjugaison pr\`es.

On \'ecrit $\pi$ comme une induite de s\'erie discr\`ete; donc le param\`etre de $\pi$ est la somme directe des param\`etres de chacune de ces s\'eries discr\`etes. On  applique la remarque pr\'ec\'edente \`a chaque s\'erie discr\`ete d'o\`u une matrice $\alpha$ de $GL(n,{\mathbb C})$ obtenue en prenant la somme directe de celles convenant pour chaque s\'erie discr\`ete. Ici on suppose uniquement que le signe associ\'e \`a chaque s\'erie discr\`ete ne d\'epend que de $\pi$ et on le note $\lambda_{\pi}$.  D'o\`u
$$
\forall w\in W_{E}\times SL(2,{\mathbb C}), \psi_{\pi,F}(w )= (\psi_{\pi,E}(w),\alpha^{-1}\, ^t\psi_{\pi,E}(w) \alpha),
$$
$$
\psi_{\pi,F}(\sigma)=\sigma_{0}(1, \lambda_{\pi}^t\alpha^{-1}\alpha).
$$
On calcule le conjugu\'e
$\psi'_{\pi,F}:=
(1,J \alpha)\psi_{\pi,F}(1,\alpha^{-1}J^{-1})$ et on a:
$$
\forall w\in W_{E}\times SL(2,{\mathbb C}), \psi_{\pi,F}(w )=(\psi_{\pi,E}(w), J ^t\psi_{\pi,E}(w)^{-1}J^{-1}); \eqno(2)
$$
$$
\psi'_{\pi,F}(\sigma)=\sigma_{0}(J\alpha,1)(\psi_{\pi,E}(\sigma^2),1)(1,\alpha^{-1}J^{-1})=
$$
$$
\sigma_{0}(J\alpha ,\lambda_{\pi} ^t \alpha^{-1}J^{-1}). \eqno(3)
$$
On remarque que les \'el\'ements de (2) sont bien invariants sous $\theta^*$. Pour (3), on remarque d'abord que $$J\, ^t(J\alpha)^{-1}J^{-1}=( J\, ^tJ^{-1} )^t \alpha^{-1}J^{-1}= (-1)^{n-1}(^t  \alpha^{-1}J^{-1}).$$ Ainsi (3) est dans le commutant de $\theta^*$ si et seulement si $\lambda_{\pi}=(-1)^{n-1}$.

Pour la r\'eciproque, on proc\`ede en sens inverse. Cela termine la preuve.
\subsection {Lien avec les fonctions $L$ d'Asai\label{asai}}
On rappelle la definition de la fonction $L$ d'Asai (cf. \cite{henniartlomeli}): un peu plus g\'en\'e\-ralement on fixe $d$ un entier et un signe $\eta$ et on consid\`ere la repr\'esentation $Asai_{\eta}$ de $(GL(d,{\mathbb C})\times GL(d,{\mathbb C}))\rtimes {\rm Gal}(E/F)$ dans ${\mathbb C}^{d^2}$ d\'efinie par:
$$\forall A\in End({\mathbb C}^d), \forall g,g' \in GL(d,{\mathbb C})
Asai_{\eta}(g,g')(A)=gA\, ^tg'
$$
et pour $\sigma_{0}$ l'\'el\'ement non trivial de ${\rm Gal}(E/F)$ et $A$ comme ci-dessus
$$
Asai_{\eta}(\sigma_{0})A=\eta\,  ^tA.
$$
\begin{cor}Soit $\pi$ une repr\'esentation $\theta$-elliptique de $GL(n,E)$; on \'ecrit $\pi=\times _{(\rho,a)\in {\mathcal E}}St(\rho,a)$ o\`u ${\mathcal E}$ est un ensemble de couples form\'es d'une  repr\'esentation cuspidale $\rho$ de $GL(d_{\rho},E)$ et d'un entier $a$. Le param\`etre de $\pi$ se factorise (\`a conjugaison pr\`es) par le $L$-groupe de $U(n,E/F)$ si et seulement si pour chaque couple $(\rho,a)\in {\mathcal E}$, la fonction $L(\sigma, \rho_{\eta},s)$ a un p\^ole en $s=0$, pour $\eta=(-1)^{n-1}$ si $a$ est impair et $\eta=(-1)^n$ si $a$ est pair.
\end{cor}
Pour ce corollaire on utilise le fait que la fonction $L$ d'Asai d\'efinie par Shahidi et par la classification de Langlands sont les m\^emes d'apr\`es \cite{henniartasai}. Ici on pourrait formuler le lemme en restant uniquement du c\^ot\'e galoisien mais le corollaire nous servira en utilisant l'\'egalit\'e des fonctions $L$.

Etant donn\'e la proposition de \ref{casunitaire} d\'ej\`a montr\'ee, il suffit de consid\'erer le  cas o\`u $\pi$ est une s\'erie discr\`ete.

Pour d\'emontrer le corollaire, on \'ecrit $\pi$ sous la forme $St(\rho,a)$ et on montre que $\lambda_{\pi}=(-1)^{a-1}\lambda_{\rho}$. On garde les notations de \ref{casunitaire}; on \'ecrit le param\`etre de $\rho$:
$$
\psi_{F,\rho}: W_{F} \rightarrow (GL(d_{\rho},{\mathbb C})\times GL(d_{\rho},{\mathbb C})\times {\rm Gal}(E/F),
$$ associ\'e au morphisme $\psi_{E,\rho}$ donn\'e par la correspondance de Langlands.

La repr\'esentation de dimension $a$ de $SL(2,{\mathbb C})$ donne lieu \`a un morphisme $\psi_{a}: SL(2,{\mathbb C}) \rightarrow GL(a,{\mathbb C})$. On obtient le morphisme associ\'e \`a $\pi$, $\psi_{E,\pi}$ en faisant le produit tensoriel, $\psi_{E,\rho}\otimes \psi_{a}$ suivi de l'inclusion de $GL(d_{\rho},{\mathbb C})\times GL(a,{\mathbb C})$ dans $GL(ad_{\rho},{\mathbb C})$. Ainsi on fixe un \'el\'ement, $\alpha_{1}$ de $GL(d_{\rho},{\mathbb C})$ qui conjugue $\psi_{E,\rho}$ et son transconjugu\'e et \'el\'ement $\alpha_{2}$ de $GL(a,{\mathbb C})$ qui conjugue $\psi_{a}$ et $^t\psi_{a}^{-1}$. L'\'el\'ement not\'e $\alpha$ dans la remarque de \ref{casunitaire} peut \^etre choisi comme valant $\alpha_{1}\times \alpha_{2}$ vu comme un \'el\'ement de $GL(d_{\rho},{\mathbb C})\times GL(a,{\mathbb C}) \hookrightarrow GL(ad_{\rho},{\mathbb C})$.

On a donc $\psi_{E,\pi}(\sigma_{0}^2)=\lambda_{\pi} \,  ^t\alpha^{-1}\alpha$ et $$^t\alpha^{-1}\alpha= \, ^t\alpha_{1}^{-1}\alpha_{1} \times \, ^t\alpha_{2}^{-1}\alpha_{2}.$$
Mais $\psi_{a}$ est orthogonal si $a$ est impair et symplectique sinon, et donc $^t\alpha_{2}^{-1}\alpha_{2}=(-1)^{a-1}$. De plus par d\'efinition de $\lambda_{\rho}$, on a
$$
\psi_{E,\rho}(\sigma_{0}^2)=\lambda_{\rho }\, ^t\alpha_{1}^{-1}\alpha_{1}.$$
Comme $\psi_{E,\pi}(\sigma_{0}^2)=\psi_{E,\rho}(\sigma_{0}^2)$, on en d\'eduit:
$$
\lambda_{\pi} \, ^t\alpha_{1}^{-1}\alpha_{1} (-1)^{a-1}=\lambda_{\rho}\, ^t\alpha_{1}^{-1}\alpha_{1}.
$$
D'o\`u $\lambda_{\pi}=\lambda_{\rho}(-1)^{a-1}$ comme cherch\'e.

Il reste \`a montrer que pour $\rho$ une repr\'esentation cuspidale autocontragr\'e\-diente, $L(\rho,Asai_{\eta},s)$ a un p\^ole en $s=0$ exactement quand $\eta=\lambda_{\rho}$. On a d\'ecrit le morphisme de $W_{F}$ dans $(GL(n,{\mathbb C})\times GL(n,{\mathbb C}))\rtimes {\rm Gal}(E/F)$ associ\'e \`a $\rho$ dans la remarque de  \ref{casunitaire}. On reprend les notations de loc. cite en rempla\c{c}ant $\pi$ par $\rho$ et on montre $Asai_{\lambda_{\rho}}\circ \psi_{\rho,F}$ laisse la matrice $^t\alpha^{-1}$ invariante. Pour cela on rappelle que pour tout $w\in W_{E}$
$$
\psi_{\rho,F}(w)=(\psi_{\rho,E}(w), \alpha^{-1}\, ^t \psi_{\rho,E}(w)^{-1}\alpha)
$$
et $Asai_{\lambda_{\rho}}(\psi_{\rho,F}(w))\, (^t\alpha^{-1})=^t\alpha^{-1}$ de fa\c{c}on \'evidente. On a pour $\sigma\in W_{F}-W_{E}$ fix\'e comme en loc. cite (une diff\'erence de choix fait varier $\alpha$), on a
$$
Asai_{\lambda_{\rho}}(\psi_{\rho,F}(\sigma))\, (^t\alpha^{-1})=Asai_{\lambda_{\rho}}(\sigma_{0})(^t\alpha^{-1}\, ^t(\psi_{\rho,E})(\sigma^2))=$$
$$
\lambda_{\\rho}Asai_{\lambda_{\rho}}(\sigma_{0})(\alpha^{-1})=^t\lambda^{-1}.
$$
Ainsi $L(\rho,Asai_{\lambda_{\rho}},s)$ a un p\^ole en $s=0$. De plus $$L(\rho\times \theta(\rho),s)=L(\rho,Asai_{+},s)L(\rho,Asai_{-},s)$$ et le premier terme n'a qu'un p\^ole d'ordre 1 en $s=0$. Ainsi $L(\rho,Asai_{-\lambda_{\rho}},s)$ n'a pas de p\^ole en $s=0$. Cela termine la preuve.

\subsection {Groupes endoscopiques  de $GL(n,E).\theta$ et image des pa\-ra\-m\`etres de Langlands autoduaux \label{morphismeendoscopique}}
L'inclusion du $L$-groupe d'un groupe endoscopique de $U(n,E/F)$ dans celui de $U(n,E/F)$ tient compte des diff\'erences de parit\'e (cf. par exemple \cite{inventaire}) et, pour \'etudier les groupes unitaires, il faut donc  choisir un caract\`ere de $E^*$ dont la restriction \`a $F^*$ est le caract\`ere quadratique correspondant \`a l'extension $E$ de $F$; on note $\omega$ un tel caract\`ere. 
Ce caract\`ere $\omega$ sert aussi \`a d\'ecrire l'inclusion des $L$-groupes des groupes endoscopiques de $GL(n,E).\theta$ dans celui de $GL(n,E)$ (cf. aussi \cite{inventaire}) et il  n'intervient pas dans les \'enonc\'es quand on les formule correctement.

La proposition ci-dessous revient au cas g\'en\'eral consid\'er\'e dans cet article.

\begin{prop}Soit $\pi$ une repr\'esentation $\theta$-elliptique de $\tilde{GL}$.

(i) Il existe un unique sous-groupe endoscopique elliptique de $\tilde{GL}.\theta$ tel que le param\`etre de $\pi$ se factorise par l'image de son $L$-groupe.

(ii)On suppose que l'unique groupe endoscopique elliptique d\'efini en (i) est simple. Alors $\pi$ n'est pas une s\'erie discr\`ete si et seulement si le param\`etre de $\pi$ se factorise par l'image du $L$-groupe d'un sous-groupe endoscopique propre de ce groupe. 
\end{prop}
Cette proposition est de l'alg\`ebre lin\'eaire dans le $L$-groupe; elle ne pose aucun probl\`eme pour les groupes orthogonaux ou symplectiques ni pour les groupes de similitudes (cf \cite{gsp4} paragraphe 2). On fait ici le cas des groupes unitaires \`a cause des signes.

On suppose donc que $\tilde{GL}=GL(n,E)$ avec $E/F$ une extension quadratique. 
On consid\`ere  les $L$ groupes comme des extensions par $W_{F}$ et pas seulement par ${\rm Gal}(E/F)$ mais les formules d\'ej\`a donn\'es par exemple pour $\psi_{\pi,E}$, $\psi_{\pi,F}$ s'\'etendent sans probl\`eme. Dans ce lemme $\pi$ n'intervient que via son param\`etre.

Supposons d'abord que $\psi_{\pi,F}$ se factorise par l'image du $L$-groupe d'un groupe endoscopique elliptique de $GL(n,E).\theta$. On a donc une d\'ecomposition $n=n_{1}+n_{2}$, le couple $(n_{1},n_{2})$ est ordonn\'e; et $\psi_{\pi,E}$ est la somme directe (ou produit) $\times_{i=1,2}\psi_{i,E}$ o\`u $\psi_{i,E}$ est \`a valeurs dans $GL(n_{i},{\mathbb C})$. Pour $i=1,2$, \`a $\psi_{i,E}$ on associe $\psi'_{i,E}$ qui est le produit tensoriel de $\psi_{i,E}$ avec un caract\`ere de $W_{E}$ correspondant si $i=1$ \`a $\omega^{n-n_{1}}$ et si $i=2$ \`a $\omega^{n-n_{2}+1}$. Et $\psi'_{i,E}$ devient un morphisme de $W_{F}\times SL(2,{\mathbb C})$ dans le $L$-groupe de $GL(n_{i},E)$ vu comme groupe sur $F$ qui se factorise par le $L$-groupe de $U(n_{i},E/F)$. Les signes associ\'es aux composantes de $\psi'_{i,E}$ sont donc $(-1)^{n_{i}-1}$. Quand on revient \`a $\psi_{i,E}$, on voit que si $i=1$ les signes sont $(-1)^{n-1}$ tandis que si $i=2$ ils valent $(-1)^{n}$: en effet avec les notations de \ref{casunitaire}, les $\alpha$ ne changent pas mais $\psi_{i,E}(\sigma^2)=\psi'_{i,E}(\sigma^2)\omega^\delta(\sigma^2)$, o\`u $\delta$ est convenable et $\omega(\sigma^2)=-1$ par d\'efinition. 

On peut renverser les arguments gr\^ace \`a \ref{casunitaire} pour montrer que $\psi_{\pi,F}$ se factorise effectivement par le groupe endoscopique d\'etermin\'e par les signes des composantes irr\'eductible de $\psi_{\pi,E}$. Et cela d\'emontre (i).

(ii):  l'hypoth\`ese de cette partie du lemme est exactement que les signes des composants de $\pi$ sont tous \'egaux; le fait que $\pi$ soit une s\'erie discr\`ete revient exactement \`a dire que $\psi_{\pi,E}$ est irr\'eductible. Si $\psi_{\pi,E}$ est irr\'eductible, il n'y a pas de factorisation par l'image d'un $L$ groupe d'un groupe endoscopique propre. Par contre si $\psi_{\pi,E}$ n'est pas irr\'eductible, n'importe quelle factorisation $\psi_{\pi,E}=\psi_{1,E}\oplus \psi_{2,E}$ donne lieu \`a une factorisation par l'image du $L$-groupe d'un groupe endoscopique propre: la d\'emonstration de (i) explique comment les signes des composants se transforment: par exemple dans le cas non principal o\`u les signes valent $(-1)^n$ ils sont multipli\'es par  en $(-1)^{n-n_{i}+1}$ et deviennent donc $(-1)^{n_{i}-1}$ pour $i=1,2$.

\section{Morphisme de Langlands des paquets stables de s\'eries discr\`etes}
\subsection{''Doubling method''\label{doublingmethod}}
La doubling method a \'et\'e initi\'ee dans le lectures notes de Gelbart, Piatetski-Shapiro et Rallis en particulier dans la premi\`ere partie de ce LN r\'edig\'ee par Cogdell \cite{GPSR}. Depuis elle a \'et\'e utilis\'ee \`a de nombreuses reprise et dans des situations tr\`es vari\'ees. Elle est d\'ej\`a utilis\'ee  avec ce m\^eme point de vue qu'ici dans \cite{book}  preuve de  6.8.1. 
\subsubsection{Le cas des groupes unitaires\label{signeunitaire}}
On commence par le cas des groupes unitaires.
\begin{lem} Soit $\pi$ une repr\'esentation cuspidale de $GL(n,E)$ autocontragr\'ediente. On suppose que le param\`etre de $\pi$ se factorise par le $L$-groupe du groupe endoscopique principal (resp. anti-principal) de ${GL}(n,E).\theta$  alors le param\`etre de $\pi$ se prolonge en un morphisme de $W_{F}\times SL(2,{\mathbb C})$ dans le $L$-groupe du groupe endoscopique anti-principal (resp. principal) de ${GL}(2n,E).\theta$. En terme de p\^ole de fonctions $L$, l'hypoth\`ese  se produit exactement quand $L(\pi,Asai_{(-1)^{n-1}},s)$ (resp. $Asai_{(-1)^n}$) a un p\^ole en $s=0$. De plus $Asai_{(-1)^{n-1}}$ (resp. $Asai_{(-1)^n})$   est la repr\'esentation du groupe dual de $\tilde{GL}(n,E)$ dans le radical nilpotent du sous groupe parabolique dual du sous-groupe parabolique du groupe endoscopique anti-principal (resp. principal) de ${GL}(2n,E).\theta$ de sous-groupe de Levi ${GL}(n,E)$.
\end{lem}
On fixe l'espace vectoriel $V$ d\'efini sur $F$ avec une structure de $E$-espace vectoriel. On suppose que $V$ est muni d'une forme hermitienne et on pose $W:=V\oplus V'$ o\`u $V'$ est l'espace $V$ muni de la forme ''oppos\'ee''. Ainsi $W$ est naturellement un espace hermitien  quasi-d\'eploy\'e. 
Ce qui nous int\'eresse d'abord est $GL(W)$ qui double $GL(V)$. 

On pose $\theta(g):=(J\, ^t\overline{g}^{-1}J^{-1})$ pour tout $g\in GL(W)$, o\`u $J$ est la matrice antidiagonale avec des $1$ et des $-1$ qui alternent. Ainsi $\theta$ laisse invariant le sous-groupe parabolique de ${GL}(W)$ qui stabilise l'espace ${V}$. Dualement $\theta^*$ agit dans le radical unipotent du parabolique ''dual'': le sous-groupe parabolique a pour Levi le groupe isomorphe  au produit de 4 copies de $GL(n,{\mathbb C})$ produit semi-direct avec ${\rm Gal}(E/F)$; pour $\sigma_{0}$ l'\'el\'ement non trivial de ${\rm Gal}(E/F)$, on a 
$$
\sigma_{0}(g,g',g_{1},g'_{1})= (g_{1},g'_{1},g,g').$$

Le radical nilpotent de l'alg\`ebre de Lie de ce sous-groupe parabolique est isomorphe \`a $End({\mathbb C}^n)\times End({\mathbb C}^n)$ avec comme action de $\sigma_{0}$:
$$
\sigma_{0}(h,h')=(h',h).$$
De plus $\theta^*$ agit aussi, en notant  $J_{n}$ la matrice antidiagonale nxn qui a aussi des 1 et des -1 qui alternent (il n'est jamais important de savoir si l'on commence par un 1 ou un -1) par:
$$
\theta^*(g,g',g_{1},g'_{1})= (J_{n}^{-1}\, ^tg_{1}^{'-1}J,J_{n}^{-1}\, ^tg_{1}^{-1}J,J_{n}^{-1}\, ^tg^{'-1}J,J_{n}^{-1}\, ^tg^{-1}J ),$$pour tout $g,g',g_{1},g'_{1}\in GL(n,{\mathbb C})$ et
$$
\theta^*(h,h')=(-J_{n}^{-1}h'J_{n},-J_{n}^{-1}hJ_{n}),$$
pour tout $h,h'\in End({\mathbb C}^n)$.

Ainsi on identifie les \'el\'ements du sous-groupe de Levi invariant sous $\theta^*$ au groupe dual de $GL(n,E)$ par:
$$
(g,g')\in GL(n,{\mathbb C})\times GL(n,{\mathbb C}) \mapsto (g, J_{n}^{-1}\, ^tg^{'-1}J_{n},g', J_{n}^{-1}\, ^tg^{-1}J_{n}),
$$
ce qui est compatible \`a l'action de ${\rm Gal}(E/F)$. Dans cette identification, une repr\'esentation $\pi$ de $GL(n,E)$ donne lieu \`a la repr\'esentation $\pi\otimes \pi^{\check\empty}$ du sous-groupe de Levi $GL(n,E)\times GL(n,E)$.

On d\'ecompose l'action de $\theta^*$ sur le radical nilpotent de l'alg\`ebre de Lie du sous-groupe parabolique en la somme de l'espace propre pour la valeur propre $+1$ et de l'espace propre pour la valeur propre $-1$, en remarquant que l'espace propre pour la valeur propre $\lambda=\pm$ est l'espace vectoriel engendr\'e par les \'el\'ements:
$
(h,-\lambda J_{n}^{-1}\, ^th J_{n}),
$ o\`u $h$ parcourt $ {\rm End}({\mathbb C}^n)$.

Ainsi chacun de ces espaces est isomorphe \`a $End({\mathbb C}^n)$ par l'isomorphisme:
$$
h\in End({\mathbb C}^n)\mapsto (hJ_{n},(-1)^n\lambda \, ^th J_{n} )$$
le signe vient de ce que $J_{n}^{-1}=\ ^tJ_{n}=(-1)^{n-1}J_{n}$. Dans ces idenfications, l'action du groupe  dual de $GL(n,E)$ est la repr\'esentation:
$$
(g,g').h= g h \, ^tg', \qquad \sigma_{0}.h=(-1)^n\lambda \, ^th.
$$
On reprend l'astuce de \cite{endoscopiqueLfunctions}: dans cet article, Arthur a remarqu\'e, bien plus g\'en\'eralement que dans le cas qui nous pr\'eoccupe, que la d\'ecomposition en espace propre sous $\theta^*$ est en fait une d\'ecomposition de la fonction $L$ associ\'e \`a $\pi$ pour l'action du $L$ groupe de $GL(n,E)$ dans le radical nilpotent du sous-groupe parabolique; en effet on fixe $z_{0}$ un \'el\'ement du centre de $GL(n,{\mathbb C})\times GL(n,{\mathbb C})$ qui agit par $-1$ sur le radical nilpotent; on prend l'\'el\'ement $(-1,1)$. L'espace propre pour la valeur propre $+1$ est exactement le radical nilpotent du sous-groupe parabolique du groupe endoscopique associ\'e au centralisateur de  $\theta^*$ et de sous-groupe de Levi $GL(n,E)$ tandis que l'autre espace propre est le radical nilpotent de l'analogue mais pour le centralisateur de $z_{0}\theta^*$. En on obtient la d\'ecomposition:
$$
L(\pi\times \, \pi^{\check\empty},s)= L(\pi,Asai_{+},s)L(\pi,Asai_{-},s),
$$
mais plus pr\'ecis\'ement on a une interpr\'etation des fonctions $L$ du membre de droite: supposons que pour $\lambda'\in \{\pm 1\}$ la fonction $L(\pi,Asai_{\lambda'},s)$ a un p\^ole en $s=0$, alors le param\`etre de Langlands de $\pi$ \`a valeurs dans le $L$-groupe de $GL(n,E)$ est dans le commutant d'un \'el\'ement nilpotent du $L$-groupe du groupe endoscopique associ\'e au centralisateur de $\theta^*$ si $\lambda'=(-1)^n$ et de $z_{0}\theta^*$ si $\lambda'= (-1)^{n-1}$. Ainsi,  si $\pi\simeq \pi^{\check\empty}$ le param\`etre de $\pi$ se factorise par le $L$ groupe du groupe endoscopique principal si et seulement si $L(\pi,Asai_{(-1)^{n-1}},s)$ a un p\^ole en $s=0$; et on vient de voir ci-dessus que le param\`etre de $\pi$ se prolonge en un morphisme de $W_{F}\times SL(2,{\mathbb C})$ dans le groupe endoscopique principal de $GL(2n,E).\theta_{2n}$ si et seulement si $L(\pi,Asai_{(-1)^{n}},s)$ a un p\^ole en $s=0$. Cela prouve le lemme.
\subsubsection{Le cas o\`u $E=F$\label{general}}
Le lemme pr\'ec\'edent est aussi vrai pour un groupe quelconque consid\'er\'e ici \`a la diff\'erence pr\`es qu'il n'y a pas que deux groupes endoscopiques elliptiques simples, la dichotomie se fait entre le groupe orthogonal (et ses variantes) et le groupe symplectique. Soit $\pi$ une repr\'esentation cuspidale irr\'eductible de $GL(n,F)\times F^*$ que l'on suppose isomorphe \`a son image sous $\theta$.
Le lemme pr\'ec\'edent devient de fa\c{c}on ici \'evidente:

 on suppose que le param\`etre de $\pi$ se factorise par un groupe orthogonal (ou de similitudes orthogonales) alors ce param\`etre se prolonge en un param\`etre de $W_{F}\times SL(2,{\mathbb C})$ dans $GL(2n,{\mathbb C})$ qui est symplectique (resp. de similitude symplectique);

on suppose que le param\`etre de $\pi$ se factorise par le groupe symplectique, ce qui n\'ecessite que $n$ soit pair, alors ce param\`etre se prolonge en un param\`etre de $W_{F}\times SL(2,{\mathbb C})$ dans $GL(2n,{\mathbb C})$ qui se factorise par le groupe $SO(2n,{\mathbb C})$. Et une assertion analogue avec le groupe des similitudes symplectiques.

\subsection{Morphisme associ\'e \`a une s\'erie discr\`ete $\theta$-stable\label{morphismediscret}}
\begin{prop}Soit $\rho$ une repr\'esentation cuspidale irr\'eductible de ${GL}(n,E)$ isomorphe \`a $\theta(\pi)$ et soit $a\in {\mathbb N}$. La trace tordue de la repr\'esentation $St(\rho,a)$ est un transfert d'un paquet de s\'eries discr\`etes du groupe endoscopique elliptique simple de $\tilde{GL}_{an}$ si et seulement si le param\`etre de Langlands de $St(\rho,a)$ se factorise par le $L$-groupe de ce groupe endoscopique.
\end{prop}

Sous les hypoth\`eses de \cite{book}, cet \'enonc\'e est 6.8.1 de \cite{book} et sans surprise c'est la m\^eme d\'emonstration que nous reprenons.

On traite d'abord le cas o\`u $a=1$ ou $2$ et le cas des groupes unitaires.
On fixe $\rho$ une repr\'esentation cuspidale de $\tilde{GL}(d_{\rho},E)$. Par autocontragr\'edience, on sait qu'il existe $\delta=0,1$ tel que la fonction $L(\rho,Asai_{(-1)^{n-1+\delta}},s)$ a un p\^ole en $s=0$. 
Le param\`etre de $\rho$ se factorise par le $L$-groupe du groupe endoscopique principal de $\tilde{GL}.\theta$ si $\delta=0$ et par le groupe endoscopique antiprincipal sinon (cf.  le lemme \ref{signeunitaire}); avec la m\^eme r\'ef\'erence on sait que $St(\rho,2)$ de $\tilde{GL}(2n,E)$ a un param\`etre qui se factorise par le $L$-groupe du groupe endoscopique antiprincipal de $\tilde{GL}(2n,E).\theta_{2n}$ si $\delta=0$ et principal sinon. Pour $i'=p$ ou $np$ on note $G_{n,i'}$ et $G_{2n,i'}$ les groupes endoscopique principaux $(i'=p)$ ou antiprincipaux $(i'=np)$. Et on pose $i=p$ si $\delta=0$ et $i=np$ si $\delta=1$.

On consid\`ere l'induite de $\rho$ d'une part \`a $G_{2n,i'}$ pour $i'=p$ ou $np$. Ces induites v\'erifient la condition de ramification d'Harish-Chandra et l'on sait donc  que l'induite de $\rho\vert\,\vert^{s}$ \`a $G_{2n,p}$ et \`a $G_{2n,np}$ est r\'eductible en un point $s=s_{0}\in {\mathbb R}$, point qui d\'epend \'evidemment du groupe. D'apr\`es les r\'esultats d'Harish-Chandra et l'interpr\'etation qu'en a donn\'e Shahidi, ce point est $s_{0}=0$ pour le groupe $G_{2n,i'}$ o\`u $i'$ est tel que $L(\rho,r_{\theta,i},s)$ n'a pas de p\^ole en $s=0$ \cite{shahidiannals} o\`u $r_{\theta,i'}$ est la repr\'esentation d\'ecrite dans le lemme \ref{signeunitaire}. Ainsi $i'=i$.
Alors, dans $G_{2n,i}$, l'induite de $\rho$ est semi-simple de longueur 2 et la diff\'erence des deux sous-repr\'esentations est une repr\'esentation elliptique. La trace de cette repr\'esentation elliptique d\'efinit un \'el\'ement de $I_{cusp}^{G_{2n,i}}$ et a donc une projection non nulle sur l'un des $I_{cusp}^{H,st}$ de \ref{pseudocoefficient} (2) pour au moins un groupe endoscopique elliptique de $G_{2n,i}$. On remarque d'abord que l'on a d\'ej\`a prouv\'e que la projection de cette trace sur $I_{cusp}^{G_{2n,i},st}$ vaut 0 (cf. la proposition de \ref{paquetstableseriesdiscretes}). Ainsi   $H$ est un groupe endoscopique elliptique propre. Par un argument par r\'ecurrence (par exemple) on v\'erifie que $H$ est n\'ecessairement de la forme $G_{n,i'}\times G_{n,i'}$. Le probl\`eme est de d\'eterminer $i'$. On  a pos\'e $\delta=0$ si $i=p$ et $\delta=1$ si $i=np$ et on d\'efinit de la m\^eme fa\c{c}on $\delta'$ en utilisant $i'$.   L'inclusion des $L$-groupes des groupes endoscopiques des groupes unitaires est d\'ecrite par exemple en \cite{inventaire}: $G_{2n,i}$ a pour groupe endoscopique $G_{n,i'}\times G_{n,i'}$ exactement quand $n-1+ (2n-n)+\delta'\equiv 2n-1+\delta[2]$, o\`u $\delta'$ vaut $0$ si $i'=p$ et $1$ sinon, c'est-\`a-dire exactement $\delta'=\delta$, ce qui est l'assertion cherch\'ee.  Il faut encore comprendre quelle repr\'esentation (a priori virtuelle) de $G_{n,i}\times G_{n,i}$ se transf\`ere en la repr\'esentation elliptique fix\'ee de $G_{2n,i}$.

On fait la m\^eme construction en rempla\c{c}ant $\rho$ par $St(\rho,3)$; on trouve encore que l'induite de $St(\rho,3)$ \`a $G_{6n,i}$ est r\'eductible n\'ecessairement semi-simple et de longueur deux et on l'\'ecrit donc sous la forme $\pi_{1,+}\oplus \pi_{2,+}$; on  remarque que pour $j=1,2$,$Jac_{\rho\vert\,\vert}Jac_{\rho\vert\,\vert}(\pi_{j,+})$ est le double de l'un des sous-module irr\'eductible de l'induite de $\rho$ (uniquement d\'etermin\'e par $j$) en particulier est non nul mais cela montre aussi que $Jac_{\rho\vert\,\vert}Jac_{\rho\vert\,\vert}(\pi_{1,+}-\pi_{2,+})$ est un multiple de la repr\'esentation elliptique consid\'er\'ee dans le paragraphe pr\'ec\'edent . 

On \'ecrit $\pi_{1,+}- \pi_{2,+}$ comme transfert de paquet stable de s\'eries discr\`etes de groupes endoscopiques elliptiques de $G_{6n,i}$. En tenant compte de \ref{pseudocoefficient} (2), il existe au moins un groupe endoscopique n\'ecessairement produit $H_{1}\times H_{2}$ et une repr\'esentation virtuelle $\Pi_{1}\otimes \Pi_{2}$ de ce groupe, n\'ecessairement combinaison lin\'eaire de s\'eries discr\`etes (cf. la proposition de \ref{paquetstableseriesdiscretes}) dont le transfert est $\pi_{1,+}-\pi_{2,+}$. On utilise la compatibilit\'e du transfert au module de Jacquet et \ref{pasdeuxfois} qui force les non nullit\'es $Jac_{\rho\vert\,\vert}\Pi_{1}\neq 0$ et $Jac_{\rho\vert\,\vert}\Pi_{2}\neq 0$; en tenant compte de \ref{supportcuspidaletendoscopie}, cela force $H_{1}=H_{2}=G_{3n,i}$ et on v\'erifie encore que le transfert de $Jac_{\rho\vert\,\vert}\pi_{1}\otimes Jac_{\rho\vert\,\vert}\pi_{2}$ est n\'ecessairement, \`a un scalaire pr\`es, une repr\'esentation elliptique de la forme $\pi_{1}\ominus \pi_{2}$ de $G_{n,i}$ o\`u $\pi_{1}\oplus \pi_{2}$ est pr\'ecis\'ement l'induite de $\rho$. On retrouve le paquet de s\'eries discr\`etes du paragraphe pr\'ec\'edent mais maintenant on sait en plus  que $Jac_{\rho\vert\,\vert}\Pi_{1}$ est un paquet stable de s\'erie discr\`ete de $G_{n,i}$ de support cuspidal \'etendu exactement \'egal \`a $\rho$. Ce qui est l'assertion cherch\'ee.

\

Le cas de $St(\rho,2)$ est beaucoup plus facile: ici on est dans le cas o\`u (avec les notations pr\'ec\'edentes) $L(\rho,Asai_{(-1)^{n-1+\delta}},s)$ a un p\^ole en $s=0$ et o\`u l'induite de $\rho\vert\,\vert^{s_{0}}$ \`a $G_{2n,i'}$ (avec $i'\neq i$) est r\'eductible en un point $s_{0}\neq 0$; dans l'induite, il y a une s\'erie discr\`ete not\'ee $\pi$.  On note $\Pi$ la repr\'esentation de $\tilde{GL}(2n,E)$ dont la trace tordue est un transfert du paquet stable contenant $\Pi$. Comme $Jac_{\rho\vert\,\vert^{s_{0}}}\pi\neq 0$, on a vu (\ref{transfertetjac}) que $Jac^{GL}_{\rho\vert\,\vert^{s_{0}}}\Pi\neq 0$; cela force $s_{0}=1/2$ et $\Pi=St(\rho,2)$. C'est exactement ce que l'on cherchait \`a montrer.

\

On consid\`ere maintenant une s\'erie discr\`ete autoduale de ${GL}(n,E)$, $St(\rho,a)$. On note $G_{an,i}$ le groupe endoscopique elliptique de ${GL}(an,E).\theta_{an}$ qui contient un paquet stable de s\'erie discr\`ete se transf\'erant en la trace tordue de $St(\rho,a)$ ou $i=p$ ou $np$ comme dans la d\'emonstration pr\'ec\'edente. En prenant des modules de Jacquet, on v\'erifie que si $a$ est impair, la trace tordue de $\rho$ est un transfert d'un paquet de s\'eries discr\`etes de $G_{n,i}$ tandis que si $a$ est pair $St(\rho,2)$ est un transfert d'un paquet de s\'eries discr\`etes de $G_{2n,i}$.  On vient alors de montrer que dans le cas o\`u $a $ est impair  le param\`etre de $\rho$ se factorise par le $L$-groupe du groupe endoscopique elliptique $G_{n,i}$ et $i$ est d\'etermin\'e par le fait que $\lambda_{\rho}=(-1)^{n-1+\delta}$ avec $\delta=0$ si $i=p$ et $\delta=1$ si $i=np$; comme $a$ est suppos\'e impair  $(-1)^{an-1+\delta}=(-1)^{n-a+\delta}$ et le param\`etre de $St(\rho,a)$ se factorise par le groupe dual de $G_{an,i}$ (d'apr\`es la proposition de \ref{casunitaire}). Si $a$ est pair, on a  vu que le param\`etre de $St(\rho,2)$ se factorise par le groupe dual de $G_{2n,i}$, c'est-\`a-dire que $\lambda_{\rho}=(-1)^{2n-1+\delta}=(-1)^{an-1+\delta}$ et le param\`etre de $St(\rho,a)$ avec $a$ pair se factorise par le groupe dual de $G_{an,i}$. Cela termine la preuve dans le cas des groupes unitaires.

\

On consid\`ere maintenant le cas $E=F$ et o\`u $\tilde{GL}=GL(n,F)$ avec  $n$  impair. Ainsi$L(\rho\times \nu,Sym^2,s)$ a un p\^ole en $s=0$, o\`u $Sym^2\otimes r_{1}$ est la repr\'esentation de $GL(n,{\mathbb C})\times {\mathbb C}^*$ dans $GL(n(n+1)/2,{\mathbb C})$. Dans ce cas, le param\`etre de $\pi $ se factorise par le groupe des similitudes orthogonales qui est isomorphe \`a $SO(n,{\mathbb C})\times {\mathbb C}^*$. Le morphisme qui se d\'eduit de $W_{F}$ dans ${\mathbb C}^*$ donne un caract\`ere qui permet de tordre $\rho$ pour la rendre autodual, c'est \'evidemment le caract\`ere ${\omega}_{\rho}\nu^{-(n-1)/2}$ o\`u ${\omega}_{\rho}$ est le caract\`ere central de $\rho$. On pose $\rho':=\rho \otimes \omega_{\rho}^{-1}\nu^{(n-1)/2}$; on montre que $\rho'$ un transfert d'un groupe symplectique avec la doubling method comme ci-dessus. Et on conclut; c'est la m\^eme d\'emonstration si il n'y a pas de facteur ${\mathbb C}^*$, le caract\`ere est d'ordre 2 mais existe aussi (c'est le caract\`ere central de $\rho$).

Quand $n$ est pair, la d\'emonstration est la m\^eme que dans le cas des groupes unitaires, il y a une dichotomie entre le fait que la fonction $L(\rho\times \nu, \wedge^2\otimes r_{1},s)$ et la fonction $L(\rho\times \nu, Sym^2\otimes r_{1},s)$ ait un p\^ole en $s=0$. Dans le deuxi\`eme cas, le param\`etre de $\rho$ se factorise par le groupe de similitude symplectique; il faut alors d\'emontrer que $\rho$ est un transfert du groupe $GSpin_{2n+1}$ tandis que dans le premier cas le param\`etre de $\rho$ se factorise par le groupe non connexe $GO(2n,{\mathbb C})$; ici on doit montrer que $\rho$ est un transfert d'une des formes de $GSpin_{2n}$ (et on a un automorphisme ext\'erieur dans ce cas); la forme de $GSpin_{2n}$ qui intervient est totalement fix\'ee par le caract\`ere $$W_{F}\mapsto GO(2n,{\mathbb C})\mapsto GO(2n,{\mathbb C})/GO(2n,{\mathbb C})^0\simeq \{\pm 1\},$$ o\`u l'exposant $0$ indique la composante neutre
comme en  \cite{book} 6.8.1. Ensuite le fait de savoir si $\rho$ est un transfert du ''bon'' groupe endoscopique, c'est-\`a-dire celui de son param\`etre ou du ''mauvais'' groupe endoscopique est la doubling method d\'evelopp\'ee dans le cas des groupes unitaires.

\subsection{Morphisme de Langlands\label{morphismedelanglands}}
\begin{thm} 
La trace tordue d'une repr\'esentation $\theta$-elliptique de $\tilde{GL}$ est un transfert d'un paquet stable de s\'eries discr\`etes d'un groupe endoscopique elliptique de $\tilde{GL}.\theta$ si et seulement si le param\`etre de la repr\'esentation de $\tilde{GL}$ se factorise par l'image du groupe dual de ce groupe endoscopique.
\end{thm}
\begin{thm}Les  paquets stables de s\'erie discr\`ete du groupe $G$  sont exactement classifi\'ees par les morphismes de Langlands de $W_{F}\times SL(2,{\mathbb C})$ dans le $L$-groupe de $G$ dont le centralisateur est un groupe fini; si $G=SO(2n,F)$ ou $GSpin_{2n}(F)$ ce sont les orbites sous l'automorphisme ext\'erieur dont nous montrons qu'elles satisfont cette bijection.
\end{thm}

Le deuxi\`eme th\'eor\`eme r\'esulte du premier. On d\'emontre donc le premier. La strat\'egie est la suivante: on fixe une repr\'esentation, $\pi^{GL}$ $\theta$-elliptique de $\tilde{GL}$ et on fixe le groupe endoscopique elliptique $\underline{H}$ par lequel le param\`etre de $\pi^{GL}$ se factorise (cf. \ref{morphismeendoscopique}). Et on va montrer que ce groupe $H$ a au moins une  repr\'esentation elliptique dont le support cuspidal \'etendu est, via l'application fix\'ee par $\underline{H}$, le support cuspidal de $\pi^{GL}$; cela suffira en vertu de \ref{deflang}. On distingue trois cas:

1e cas: on suppose que $H$ est un produit de groupes; on applique alors le deuxi\`eme th\'eor\`eme par r\'ecurrence.

2e cas: on suppose que $H$ est un groupe simple. Si $\pi^{GL}$ est une s\'erie discr\`ete, on a d\'ej\`a d\'emontr\'e le th\'eor\`eme. Sinon, on consid\`ere la d\'ecomposition du param\`etre de $\pi^{GL}$ vu comme repr\'esentation de $W_{F}\times SL(2,{\mathbb C})$ en sous-repr\'esentation irr\'eductible; on factorise cette d\'ecomposition en somme de 2 sous-repr\'esentations donnant lieu \`a une factorisation de ce param\`etre par un sous-groupe de $H$ , \`a une exception pr\`es (qui fera l'objet du 3e cas) si cette d\'ecomposition est la somme de deux repr\'esentations orthogonales (ou de similitude orthogonales) de dimension impaire;  ce sous-groupe fait partie d'une donn\'ee endoscopique elliptique $\underline{H}'$ de $H$.il existe alors au moins une donn\'ee endoscopique elliptique $\underline{H}'$ de $\underline{H}$. On applique l'hypoth\`ese de r\'ecurrence \`a $\underline{H}'$, d'o\`u un paquet stable de s\'erie discr\`ete de $\underline{H}'$ dont le support cuspidal \'etendu est par construction le support cuspidal de $\pi^{GL}$. On transf\`ere ce paquet de s\'eries discr\`etes \`a $\underline{G}$ et l'image est une combinaison lin\'eaire de s\'eries elliptiques qui ont pour support cuspidal \'etendu l'image du support cuspidal \'etendu du paquet stable de s\'erie discr\`etes de $\underline{H}'$ d'apr\`es \ref{supportcuspidaletendoscopie} et c'est donc bien des repr\'esentations elliptiques ayant le support cuspidal cherch\'e.

3e cas: on suppose que $H$ est le groupe $SO(2n,F)$ ou $GSpin(2n,F)$ et que le param\`etre de $\pi^{GL}$ est la somme de deux repr\'esentations irr\'eductibles chacune de dimension impaire. En prenant des modules de Jacquet, on se ram\`ene au cas o\`u $\pi^{GL}= (\rho\times \rho')\times \nu,
$
o\`u $\rho$ et $\rho'$ sont des repr\'esentations cuspidales unitaires de $GL(d_{\rho},F)$ et $GL(d_{\rho'},F)$ avec $d_{\rho}+d_{\rho'}=2n $ et $d_{\rho},d_{\rho'}$ des entiers impairs. Pour traiter ce cas, on reprend la preuve de \ref{morphismediscret}: on consid\`ere l'induite de la repr\'esentation $(\rho\times \rho')\times \nu$ de $GL(2n,F)\times F^*$ \`a chacun des groupes endoscopiques elliptiques simples de $GL(4n,F)\times F^*$. L'induite \`a $GSpin(4n,F)$ est r\'eductible par les r\'esultats d'Harish-Chandra; en effet dans ce cas, le $R$-groupe est n\'ecessairement non trivial (cf. pour le cas de $SO(4n,F)$, le travail de Goldberg \cite{goldbergclassique} d'o\`u vient cette remarque) car le $R$-groupe est le quotient du sous-groupe du groupe de Weyl stabilisant la repr\'esentation par un sous-groupe n\'ecessairement engendr\'e par des sym\'etries \'el\'ementaires. Mais la condition de parit\'e entra\^{\i}ne qu'aucune sym\'etrie \'el\'ementaire ne stabilise la repr\'esentation que l'on induit. L'induite est alors de dimension 2 et la diff\'erence des deux sous-modules irr\'eductibles est une repr\'esentation elliptique. Ensuite on conclut comme dans la preuve de \ref{morphismediscret} pour montrer que la trace de cette repr\'esentation est un transfert d'un paquet de s\'eries discr\`etes de $GSpin(2n,F)$ qui ont le bon support cuspidal \'etendu.
\section{$R$-groupe et cardinal des paquets stables de repr\'esentations temp\'er\'ees\label{cardinal}}
\subsection{Lien des constructions pr\'ec\'edentes avec les R-groupes}
La d\'etermination des $R$-groupes a fait l'objet de nombreux articles, en particulier leur d\'etermination explicite avec les travaux de R. Herb et D. Goldberg \cite{goldberg}. Dans le th\'eo\`eme et la remarque suivante, on inclut bien le cas des groupes orthogonaux pairs et de $GSpin^*_{2n}$, c'est-\`a-dire les groupes non connexes. 
\begin{thm}(\cite{book}, pour certains groupes)
Soit $\pi$ une s\'erie discr\`ete irr\'eductible de $G$ et $Jord(\pi)$ son ensemble de bloc de Jordan.

(i) Le nombre d'\'el\'ement du paquet de Langlands contenant $\pi$ est exactement $2^{\vert Jord(\pi)\vert-1}$.

(ii) Soit $\rho$ une repr\'esentation irr\'eductible cuspidale et unitaire de $GL(d_{\rho},E)$  et soit $a$ un entier.  L'induite $St(\rho,a)\times \pi$ est r\'eductible, n\'ecessairement de longueur deux, si et seulement si 

$\rho\simeq \theta(\rho)$, $(\rho,a)\notin Jord(\pi)$ et  $L(\rho,r_{G},s)$ a un p\^ole en $s=0$ si $a$ est pair et n'en a pas si $a$ est impair.
\end{thm}

Montrons le th\'eor\`eme. 

Pour (i), on a montr\'e que l'endoscopie respecte le support cuspidal \'etendu. Pour $\underline{H}$ une donn\'ee endoscopique elliptique de la composante neutre de $G$, on note $I_{cusp, Jord(\pi)}^{\underline{H},st}$ l'espace des combinaisons lin\'eaires stables bas\'ees form\'es de repr\'esentation ayant pour support cuspidal \'etendu d'image $Jord(\pi)$ dans ${GL}(n,{\mathbb C})$. Et on note $I^0_{cusp,Jord(\pi)}$ le sous-espace vectoriel de $I_{cusp}$ engendr\'e par les traces des s\'eries discr\`etes dans le paquet de Langlands de $\pi$; pour le cas de $O(2n,F)$ et $GSpin_{2n}^*(F)$, on consid\`ere  le groupe connexe, $SO(2n,F)$ et $GSpin_{2n}(F)$. On a donc
$$
dim\, I^0_{cusp,Jord(\pi)}=\sum_{\underline{H}/\simeq}dim\, I_{cusp,Jord(\pi)}^{\underline{H},st,{Aut(\underline{H})}},
$$avec les notations de \ref{pseudocoefficient}.

On a vu que $dim\, I_{cusp,Jord(\pi)}^{G,st, Aut(G)}=1$ (\ref{deflang}). On fixe $\underline{H}$ une donn\'ee endoscopique elliptique propre et on utilise la description de \cite{inventaire}: le groupe sous-jacent est le produit de deux groupes du type de ceux que l'on \'etudie ici (ou de leur composante connexe). Ils ont des paquets stables de s\'erie discr\`etes de support cuspidal \'etendu (d'image) $Jord(\pi)$ si et seulement si $Jord(\pi)$ se d\'ecoupe en deux sous-ensembles relatifs chacun \`a l'un des groupes. 

Dans le cas des groupes unitaires en dimension impair, il n'y a pas de groupes d'automorphisme d'une donn\'ee endoscopique. On termine ce cas: pour chaque d\'ecoupage de $Jord(\pi)$ en deux sous-ensembles non ordonn\'e et propre, il existe exactement un couple de donn\'ees endoscopiques, conjugu\'ee l'une de l'autre par un automorphisme ayant un paquet stable de repr\'esentation avec ce support cuspidal \'etendu. En ajoutant le couple form\'e de $G$ lui-m\^eme et de l'ensemble vide, on voit que le cardinal cherch\'e est exactement le nombre de d\'ecoupage de $Jord(\pi)$ en deux sous-ensemble non ordonn\'e mais dont l'un peut \^etre l'ensemble vide. Cela est exactement $2^{\vert Jord(\pi)\vert-1}$.

Dans le cas des groupes unitaires en dimension paire et des groupes or\-thogonaux ou de similitude $GSpin$ en dimension impaire,
il n'y a d'automorphis\-mes que  si les deux groupes sont \'egaux (ceci ne pouvait pas se produire dans l'exemple pr\'ec\'edent).  On raisonne comme dans le cas pr\'ec\'edent, tout d\'ecoupage de $
Jord(\pi)$ en deux sous-ensemble dont aucun n'est vide, d\'etermine un couple de donn\'ees endoscopiques, ici non n\'ecessairement distincts,  qui sont conjugu\'e l'une de l'autre par automorphisme non trivial;  si les donn\'ees sont confondus, le d\'ecoupage de $Jord(\pi)$ se fait quand m\^eme en deux ensembles \'evidemment distincts et donne lieu \`a un espace de dimension un puisque le d\'ecoupage est fait \`a l'ordre pr\`es. On a donc le m\^eme calcul que ci-dessus.

Il reste le cas des groupes symplectiques trait\'e par \cite{book}, on ne revient donc pas dessus et des groupes orthogonaux ou $GSpin$ en dimension paire. Le cas des groupes orthogonaux est totalement d\^u \`a \cite{book} dont en particulier le chapitre 8 sp\'ecifique \`a ce cas. On donne un argument l\'eg\`erement diff\'erent en montrant la remarque ci-dessous qui est aussi dans \cite{book}; dans la remarque ci-dessous $G$ est le groupe non connexe $O(2n,F)$ ou $GSpin^*(2n,F)$:
\begin{rmq}
Soit $\tau$ une s\'erie discr\`ete de $G$, la restriction de $\tau$ \`a la composante neutre de $G$ est irr\'eductible sauf exactement quand $Jord(\tau)$ n'a que des couples $(\rho',a')$ tels que $a'd_{\rho'}$ soit pair.
\end{rmq}

En effet, les s\'eries discr\`etes de $GSpin^*_{2n}(F)$ et $O(2n,F)$ qui interviennent dans l'endoscopie tordue de $GSpin^*_{2n}(F)$ et $O(2n,F)$ sont exactement celles dont la restriction \`a la composante neutre reste irr\'eductible. Les groupes endos\-copiques pour cette endoscopie tordue ont un groupe dual dont la composante neutre est un produit de groupes orthogonaux impairs; cette endoscopie tordue respectent elle aussi le support cuspidal \'etendu. Ainsi toute s\'erie discr\`ete tel que $Jord(\pi)$ ne contient que des couples $(\rho',a')$ avec $a'd_{\rho'}$ pair, n'est pas la restriction d'une s\'erie discr\`ete irr\'eductible du groupe non connexe.

La preuve de la remarque sera termin\'ee ci-dessous mais on termine, dans ce cas, la preuve du th\'eor\`eme: on suppose que $Jord(\pi)$ a la propri\'et\'e de parit\'e ci-dessus et il en est donc ainsi pour toutes les repr\'esentations intervenant dans $I_{cusp, Jord(\pi)}^{\underline{H},st}$. Ainsi $I_{cusp, Jord(\pi)}^{0,G,st}$ est de dimension 2. Les automorphismes entre les donn\'ees endoscopiques inclus l'\'echange des deux facteurs et l'action du produit des automorphismes ext\'erieurs non triviaux. 
Et on trouve alors que $dim\, I_{cusp,Jord(\pi)}^0$ vaut deux fois le nombre de d\'ecoupage de l'ensemble $Jord(\pi)$ en deux sous-ensembles non ordonn\'es dont l'un peut \^etre vide, c'est-\`a-dire vaut $2^{\vert Jord(\pi)\vert}$. Quand on repasse au groupe non connexe, d'apr\`es ce que l'on vient de voir le cardinal est divis\'e par deux et on trouve l'\'enonc\'e du th\'eor\`eme.

\

On suppose maintenant que $Jord(\pi)$ contient au moins un couple $(\rho,a)$ tel que $ad_{\rho}$ est impair. On note $m_{i}$ le nombre d'\'el\'ements de $Jord(\pi)$ ayant cette propri\'et\'e de parit\'e et $m_{p}$ le nombre d'\'el\'em\'ents de $Jord(\pi)$ ayant la parit\'e oppos\'ee. On va utiliser le fait que le nombre de d\'ecoupage d'un ensemble ayant un nombre, $m_{0}$,  pair d'\'el\'em\'ents en deux sous-ensemble non ordonn\'e ayant un nombre impair d'\'el\'ements est $2^{m_{0}-1}$. Alors
l'endoscopie pour le groupe non connexe montre que
$
dim\, I^0_{cusp,Jord(\pi)}=2x+2^{m_{p}}2^{m_{i}-1}/2$, o\`u $x$ est le nombre de s\'eries discr\`etes de $G$, associ\'e \`a $Jord(\pi)$,  dont la restriction \`a $G^0$ se coupe en deux et o\`u la division par $2$ vient des automorphismes (comme plus haut) pour l'endoscopie de $G$.

On admet le r\'esultat par r\'ecurrence sur le rang de $G$ pour les groupes endoscopiques propres de la composante neutre de $G$.  Ainsi pour une telle donn\'ee endoscopique propre,  l'action du produit des automorphismes ext\'erieurs non triviaux sur $I_{cusp,Jord(\pi)}^{\underline{H},st}$  est  trivial sur au moins un des facteurs, le support cuspidal cuspidal contenant au moins un  $(\rho',a')$ avec $a'd_{\rho'}$ impair et que la dimension des invariants sous ce groupe est exactement le nombre de d\'ecoupage de $Jord(\pi)$ adapt\'e \`a $\underline{H}$. On obtient donc 
$$
dim\, I^0_{cusp,Jord(\pi)}=y+2^{m_{p}}2^{m_{i}-1}/2,$$
o\`u $y$ vaut $0$ si la remarque est vraie, c'est-\`a-dire si la distribution stable associ\'ee \`a $G$ ne se coupe pas en restriction \`a $G^0$ et $1$ sinon. 
En comparant, on trouve $y=2x$ et par parit\'e cela force, $y=x=0$, d'o\`u la remarque. Et on trouve ici que le nombre de s\'eries discr\`etes de $G^0$ avec le support cuspidal fix\'e est $2^{\vert Jord(\pi)\vert-2}$; et il y a \'evidemment deux s\'eries discr\`etes de $G$ ayant m\^eme restriction \`a $G^0$, elles se d\'eduisent par un caract\`ere signe d'o\`u le (i) du th\'eor\`eme.

\

Montrons le (ii) du th\'eor\`eme: la condition $\rho\simeq \theta(\rho)$ est n\'ecessaire pour avoir r\'eductibilit\'e: dans le cas o\`u $G$ est connexe, c'est la condition d'Harish-Chandra, dans le cas o\`u $G$ n'est pas connexe, c'est soit la condition d'Harish-Chandra, soit la condition de Mackey; on suppose donc dans la suite que $\rho\simeq \theta(\rho)$.

On fixe $(\rho,a)$ comme dans l'\'enonc\'e mais on ne fixe pas $\pi$, on fixe uniquement $Jord(\pi)$ que l'on note ${\mathcal E}$. On consid\`ere l'ensemble des repr\'esentations induites $St(\rho,a)\times \pi'$, o\`u $\pi'$ est une s\'erie discr\`ete telle que $Jord(\pi')={\mathcal E}$. Cette induite est r\'eductible si et seulement si elle a deux sous-modules irr\'eductible dont la diff\'erence est une repr\'esentation elliptique; il faut donc ici compter le nombre de repr\'esentations elliptiques dont le support cuspidal \'etendu est ${\mathcal E}\cup\{(\rho,a),(\rho,a)\}$. Si $G$ est connexe,  il faut sommer $dim\, I_{cusp,{\mathcal E}\cup\{(\rho,a),(\rho,a)\}}^{\underline{H},st}$  pour toute donn\'ee endoscopique elliptique   $\underline{H}$ de $G$ prise \`a automorphisme pr\`es  
Pour cela il faut que dans l'ensemble ${\mathcal E}\cup\{(\rho,a),(\rho,a)\}$, $(\rho,a)$ qui intervient avec multiplicit\'e au moins 2 n'intervienne pas avec multiplicit\'e 3; d'o\`u la n\'ecessit\'e de la condition $(\rho,a)\notin Jord(\pi)$; on trouve aussi la n\'ecessit\'e de la derni\`ere condition. Si ces conditions sont satisfaites, on v\'erifie comme ci-dessus, que le nombre de repr\'esentations elliptiques cherch\'ees est exactement le nombre de d\'ecoupage de ${\mathcal E}$ en deux sous-ensembles non ordonn\'ees. Ainsi les conditions sont aussi suffisantes.

On consid\`ere maintenant le cas o\`u $G$ est non connexe: si $ad_{\rho}$ est pair, la d\'emonstration est analogue \`a celle que l'on vient de faire. On suppose donc que $ad_{\rho}$ est impair.

On regarde d'abord le cas o\`u la restriction de $\pi$ au groupe connexe et la somme $\pi_{1}\oplus \pi_{2}$ de deux repr\'esentation. D'apr\`es ce que l'on a vu ci-dessus, c'est exactement le cas quand tout \'el\'ement de $Jord(\pi)$ est form\'e de termes $(\rho',a')$ avec $a'd_{\rho'}$ pair. Ainsi toutes les conditions de l'\'enonc\'e sont automatiquement satisfaites et il faut donc d\'emontrer que la repr\'esentation induite $St(\rho,a)\times \pi$ est r\'eductible.

 D'apr\`es les r\'esultats g\'en\'eraux sur le $R$-groupe d'Harish-Chandra, les induites, pour $i=1,2$, $St(\rho,a)\times \pi_{i}$ sont irr\'eductibles mais elles sont conjugu\'ees l'une de l'autre par un \'el\'ement de $G_{n+2ad_{\rho}}^0$, la composante neutre; elles sont donc isomorphes. Ainsi  l'induite de $St(\rho,a)\times \pi_{1}$ \`a $G_{n+2ad_{\rho}}$ est $St(\rho,a)\times \pi$ et en restriction \`a la composante neutre est la somme de deux repr\'esentations isomorphes. Cela prouve que $St(\rho,a)\times \pi$ est la somme de deux repr\'esentations irr\'eductibles qui diff\`erent par le caract\`ere signe; on a bien montr\'e la r\'eductibilit\'e. 
 
 On consid\`ere maintenant le cas o\`u la restriction de $\pi$ au groupe connexe est irr\'eductible. On note encore ${\mathcal E}:=Jord(\pi)$. On compte le nombre de repr\'esen\-tations elliptiques de la composante neutre de $G$ ayant pour support cuspidal \'etendu ${\mathcal E}\cup\{(\rho,a),(\rho,a)\}$; on trouve que cela vaut exactement le nombre de d\'ecoupage de ${\mathcal E}$ en deux sous-ensemble non ordonn\'e et dont chacun contient un nombre impair de termes  $(\rho',a')$ avec $a'd_{\rho'}$ impair. C'est donc exactement le nombre de s\'eries discr\`ete de la composante neutre de $G$ ayant comme support cuspidal \'etendu ${\mathcal E}$. Ainsi pour toute s\'erie discr\`ete $\pi_{0}$ de la composante neutre de $G$, l'induite $St(\rho,a)\times \pi_{0}$ est r\'eductible; l'induite de cette repr\'esentation \`a $G$ n'est autre que $St(\rho,a)\times \pi$ qui est donc aussi n\'ecessairement r\'eductible. Cela termine la d\'emonstration.
\subsection{Unicit\'e de la d\'efinition des blocs de Jordan \label{unicitedefinition}}
\begin{rmq} Le th\'eor\`eme pr\'ec\'edent montre que les blocs de Jordan tel que d\'efinis ici sont bien ceux de \cite{europe}, \cite{jams}, \cite{pourshahidi}. On a aussi la caract\'erisation suivante sugg\'er\'ee par M. Tadic: $(\rho,a)\in Jord(\pi)$ si et seulement si $St(\rho,a)\times \pi$ est irr\'eductible mais il existe un entier $b$ de m\^eme parit\'e que $a$ tel que la repr\'esentation induite $St(\rho,b)\times \pi$ est r\'eductible.
\end{rmq}

\end{document}